\documentclass[11pt]{article}%
\usepackage{amssymb}
\usepackage{amsfonts}
\usepackage{amsmath}
\usepackage{graphicx}
\usepackage{algorithmic}
\usepackage{float}
\usepackage[ruled,nothing]{algorithm}
\usepackage{subfig}%
\providecommand{\U}[1]{\protect\rule{.1in}{.1in}}
\newtheorem{theorem}{Theorem}

\newtheorem{criterion}[theorem]{Criterion}
\newtheorem{definition}[theorem]{Definition}

\newenvironment{proof}[1][Proof]{\noindent\textbf{#1.} }{\ \rule{0.5em}{0.5em}}
\newdimen\dummy
\dummy=\oddsidemargin
\begin{document}

\title{Reverse Line Graph Construction: The Matrix Relabeling Algorithm MARINLINGA
Versus Roussopoulos's Algorithm}
\author{D. Liu, S. Trajanovski, P. Van Mieghem\thanks{Faculty of Electrical
Engineering, Mathematics and Computer Science, Delft University of Technology,
P.O Box 5031, 2600 GA Delft, The Netherlands; \emph{email}: \{D.Liu,S.
Trajanovski,P.F.A.VanMieghem\}@tudelft.nl. }}
\date{Delft University of Technology\\
v8, \today}
\maketitle

\begin{abstract}
We propose a new algorithm MARINLINGA for reverse line graph computation,
i.e., constructing the original graph from a given line graph. Based on the
completely new and simpler principle of link relabeling and endnode
recognition, MARINLINGA does not rely on Whitney's theorem while all previous
algorithms do. MARINLINGA has a worst case complexity of $O(N^{2})$, where $N$
denotes the number of nodes of the line graph. We demonstrate that MARINLINGA
is more time-efficient compared to Roussopoulos's algorithm, which is
well-known for its efficiency.

\end{abstract}

\floatstyle{ruled} \newfloat{algorithm}{htbp}{loa}
\floatname{algorithm}{Algorithm}


\section{Introduction}

The line graph $l\left(  G\right)  $ of a graph $G$\ is a graph in which every
node corresponds to a link of $G$ and two nodes are adjacent if and only if
their corresponding links are adjacent in $G$ (two links are adjacent if they
are incident to the same node). The graph $G$ is called the original or root
graph of $l\left(  G\right)  $. There exist examples of line graphs from
social network. Given $M$ clubs and $N$ students at an university, every
student joins two clubs. Each student has different choices (we assume that
there are enough clubs). We define two networks $G_{1}$ and $G_{2}$. The $M$
clubs are the nodes of $G_{1}$ and two nodes are adjacent if two clubs have
the same student as their member. The $N$ students are the nodes of $G_{2}$
and two nodes are adjacent if two students belong to the same club Clearly,
$G_{2}$ is the line graph of $G_{1}$. Such pairs $\left(  G_{1},G_{2}\right)
$ are common in on-line social networks like Facebook, Twitter and etc., where
users join the special groups where they share the same interest with others.
Computing the line graph of a graph and constructing the original graph of a
line graph also play an important role in link partitioning of communities
\cite{PhysRevE.80.016105_linkPartition}\cite{Nature_1}\cite{Evans_2}%
\cite{Manka}, bond percolation threshold predictions
\cite{PhysRevE.75.011114_bondPercolationThresholdPredictions}, and it also
enables us to compare the properties of a random line graph \cite{RLG} and its
original graph.

The following formula \cite{PVM_SpectraCUP} can be used to compute the
adjacency matrix of the line graph $l\left(  G\right)  $ of a graph $G$,
\begin{equation}
A_{l\left(  G\right)  }=\left(  R^{T}R\right)  _{L\times L}-2I
\label{adjacency_matrix_line_graph}%
\end{equation}
where $R$ is the incidence matrix of the undirected graph $G$. If link $j$ is
incident to node $i$, the entry $r_{ij}$ of $R$ is $1$, otherwise $0$. In each
column there are exactly two $1$-entries.

\begin{figure}
[ptb]
\begin{center}
\includegraphics[
height=2.5728in,
width=3.0078in
]%
{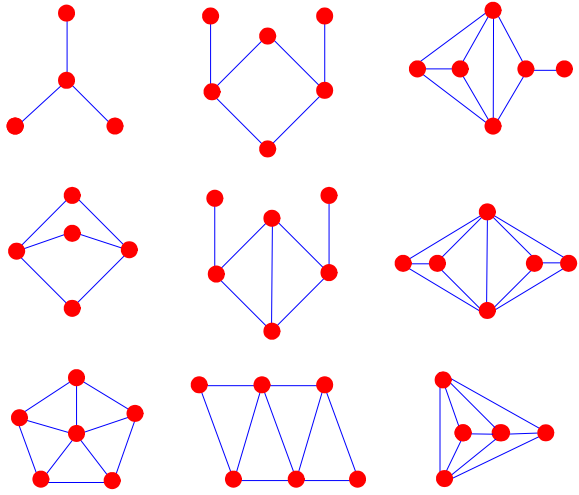}%
\caption{The nine forbidden subgraphs for line graphs \cite{Beineke}.}%
\label{ninefbdgrs}%
\end{center}
\end{figure}

Constructing the original graph is far more complex than computing the line
graph. Before constructing the original graph from a given graph, it is
important to know whether the graph is a line graph. Up till now, the
following criteria for a graph to be a line graph exist in the literature:

\begin{itemize}
\item A graph is a line graph if and only if it is possible to find a
collection of cliques in the graph, partitioning all the links, such that each
node belongs to at most two of the cliques (some of the cliques can be a
single node) and two cliques share at most one node \cite{Harary}. If the
graph is not $K_{3}$, there can be only one partition of this type.

\item A graph is a line graph if and only if it does not have the complete
bipartite graph $K_{1,3}$ as an induced subgraph, and if two odd
triangles\footnote{If every node is adjacent to two or zero nodes of a
triangle, it is an even triangle.} have a common link, the subgraph induced by
their nodes is the complete graph $K_{4}$ \cite{Rooij}.

\item A graph is a line graph if and only if none of the nine forbidden
subgraphs (see Figure \ref{ninefbdgrs}) is an induced subgraph of it
\cite{Beineke}.

\item A graph is not a line graph \cite{PVM_SpectraCUP} if the smallest
eigenvalue of the adjacency matrix $\left(  \ref{adjacency_matrix_line_graph}%
\right)  $ is smaller than $-2$.
\end{itemize}

%

The complete graph on three nodes $K_{3}$ is a line graph, which has two
different original graphs, $K_{3}$ and $K_{1,3}$ (Figure \ref{s1smaller4_2}
(b)). Except for $K_{3}$, Whitney's theorem \cite{Whitney}\cite{Harary} states
that, all line graphs have only one original graph (isomorphic graphs are
considered as the same graph). Based on the above criteria and Whitney's
theorem, several algorithms for constructing the original graph have been
proposed \cite{Lehot}\cite{Roussopoulos}\cite{Naor}\cite{Degiorgi}. Among
those algorithms, Roussopoulos's \cite{Roussopoulos} and Lehot's \cite{Lehot}
solutions are worth mentioning here.

Roussopoulos's algorithm starts with choosing an arbitrary link in the input
graph and calculating the number of triangles containing this link. Depending
on this value the starting cell is determined. The starting cell is a complete
graph $K_{m}$; if $m=2$ it is a link; if $m=3$ a triangle that contains the
starting link. Having a starting cell of the input graph, the algorithm of
Roussopoulos continues to find a clique, which is deleted. In addition, in
each step the vertices of the clique are labeled by a group number. One node
in a line graph cannot be assigned to more than two groups (otherwise it is
not a line graph). The nodes of the original graph are those partitions and
all nodes are assigned to exactly one partition. In the constructed graph
there is a link between two nodes, if the nodes are assigned to partitions
that have a non-empty intersection. The approach of Roussopoulos is based on
finding the largest connected components and sequentially the number of
triangles that contain this link. Theoretically finding the largest connected
component is, however, an $NP$-complete problem \cite{VISMARA}. Lehot's
solution is based on the characterization of line graphs by Van Rooij and Wilf
\cite{Lehot}\cite{Rooij}.

In this paper, we propose a new algorithm, the \textit{MA}trix \textit{R}%
elabeling \textit{IN}verse \textit{LIN}e \textit{G}raph \textit{A}lgorithm, in
short \textit{MARINLINGA}, that constructs the original graph given the line
graph. MARINLINGA does not explicitly rely on Whitney's theorem, as all
previous companion algorithms, but uses link relabeling and endnode
recognition in a new way. Via extensive simulation analysis, we have compared
MARINLINGA with Roussopoulos's algorithm. We demonstrate that MARINLINGA
consumes less CPU running time. The algorithms are tested on the same
machine\footnote{Processor Intel Core $2$ Duo CPU T$9600$ @ $2.80$ GHz and
$2.96$ GB RAM memory on Java Execution Environment JAVA-SE $1.6$\emph{ }and
Eclipse IDE (version Galileo $3.5$).} and we use the same input line graphs
for both algorithms.

\section{Link adjacency matrix (LAM) and line graph}

\label{LAM_and_line_graph}Two nodes of a graph are said to be adjacent if
there is a link directly connecting them. The adjacency matrix $A$ of a graph
contains all information of node adjacency: if node $i$ and node $j$ are
adjacent, the entry $a_{ij}=1$, otherwise $a_{ij}=0$. Similarly, two links are
adjacent if they are incident to the same node.

\begin{definition}
The link adjacency matrix (LAM) $C$ of a graph $G$ with $N_{G}$ nodes and
$L_{G}$ links is the $L_{G}\times L_{G}$ symmetric matrix with the entry
$c_{ij}=1$ if link $i$ and link $j$ of $G$ are adjacent, else $c_{ij}=0$.
\end{definition}

The line graph $l\left(  G\right)  $ of the graph $G$ has $N_{l\left(
G\right)  }$ nodes and $L_{l\left(  G\right)  }$ links, and consequently we
have $L_{G}=N_{l\left(  G\right)  }$. According to the definitions of the line
graph and the LAM, evidently, the LAM $C$ of $G$ is equal to the adjacency
matrix $A_{l\left(  G\right)  }$ of $l\left(  G\right)  $,%
\begin{equation}
C=A_{l\left(  G\right)  }%
\end{equation}

Due to Whitney's theorem and ignoring isomorphisms, for any graph except
$K_{3}$ and $K_{1,3}$, one can construct the graph exclusively from its LAM.
Usually, the (node) adjacency matrix is used to represent a graph. Here we use
the LAM to specify any graph, except for $K_{3}$ and $K_{1,3}$. Constructing
the original graph of a line graph is equivalent to converting a graph
representation from the LAM to the adjacency matrix. By constructing the
original graph directly from the line graph, confusion will arise concerning
the links in the original graph and the nodes in the line graph. By
introducing the concept of LAM, we can avoid confusion and facilitate the
description of our algorithm MARINLINGA.

\section{Properties of the LAM}

For a simple (undirected, unweighted and without self-loops) graph $G\left(
N_{G},L_{G}\right)  $ with $N_{G}$ nodes and $L_{G}$ links, the LAM $C$ has
more constraints than the corresponding adjacency matrix $A$, besides being
symmetric and containing only $0$ and $1$ entries.

\begin{figure}
[ptb]
\begin{center}
\includegraphics[
height=1.2877in,
width=3.1955in
]%
{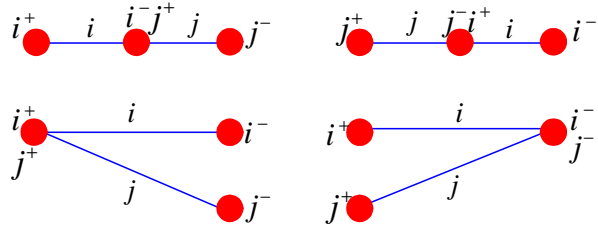}%
\caption{The four possible configurations in which link $i$ is adjacent to
link $j$.}%
\label{link_adjacent_possibility}%
\end{center}
\end{figure}

A link $i$ has two endnodes, the left endnode $i^{+}$ and the right endnode
$i^{-}$. Link $j$ also has endnodes $j^{+}$ and $j^{-}$. There are four
configurations where link $i$ is adjacent to link $j$, as shown in Figure
\ref{link_adjacent_possibility}. For each single pair of links, the LAM only
indicates whether they are adjacent. If they are adjacent, we still do not
know in which of the four possible configurations this pair of links is
adjacent. Fortunately, by combining the adjacency relation of $3$ or more
links, we can determine the configuration of those links.%

\begin{definition}
If $m$ links ($m\geq2$) are adjacent to link $i$ and incident to the same
endnode of link $i$, these $m$ links are pairwise adjacent.
\label{definition_transitive_adjacency}
\end{definition}

\begin{definition}
The links, which are adjacent to link $i$, are defined as the neighboring
links of link $i$. \label{definition_neighboring_links}
\end{definition}

\begin{definition}
The links incident to the left endnode $i^{+}$ of a link $i$ are defined as
the left-neighboring links of $i$, and the links incident to the right endnode
$i^{-}$ are defined as the right-neighboring links of $i$.
\label{definition_left&right_neighboring_links}
\end{definition}

If we can recognize the link adjacency pattern of a link and its neighboring
links, we can specify LAM entirely.

\begin{figure}
[ptb]
\begin{center}
\includegraphics[
height=2.444in,
width=4.7072in
]%
{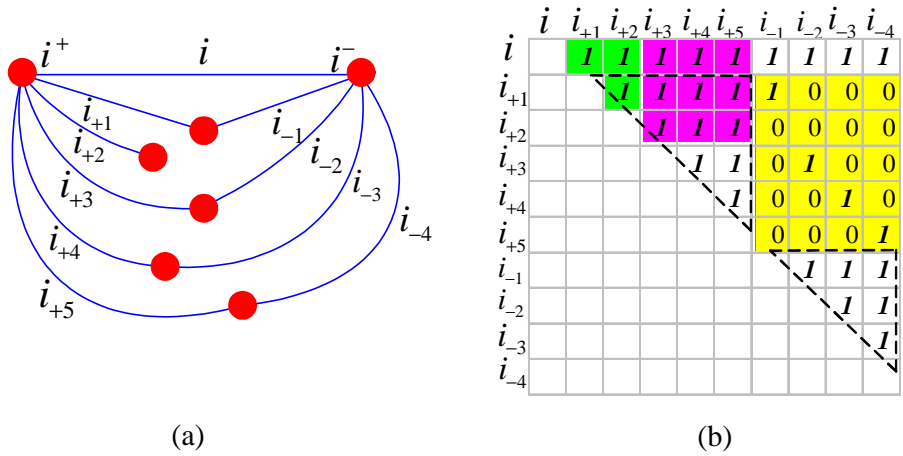}%
\caption{(a) The configuration of a link $i$ and its neighboring links. (b)
The corresponding link adjacency pattern. there is at most one $1$-entry in
each row/column of the submatrix in yellow. If all the entries in green and
magenta are $1$-entries, the entries of the triangle in white must be also
$1$-entries.}%
\label{link_adjacency_init}%
\end{center}
\end{figure}

Figure \ref{link_adjacency_init} (a) depicts an example of a link and its
neighboring links. The link $i$ has $5$ left neighboring links at its left
endnode $i^{+}$, denoted as $i_{+1},\cdots,i_{+5}$, and $4$ right neighboring
links at its right endnode $i^{-}$, denoted as $i_{-1},\cdots,i_{-4}$. The
link adjacency pattern of these $10$ links is shown in Figure
\ref{link_adjacency_init} (b). In the link adjacency pattern, the labels of
the left-neighboring links $i_{+1},\cdots,i_{+5}$ are larger than link $i$,
and smaller than the right-neighboring links $i_{-1},\cdots,i_{-4}$.

Given the configuration of link $i$ and its neighboring links, the
corresponding link adjacency pattern conforms to the following rules:

\begin{enumerate}
\item the left-neighboring links (such as $i_{+1},\cdots,i_{+5}$ in the
example of Figure \ref{link_adjacency_init} (a)) are incident to the same
endnode $i^{+}$, and are said (Definition
\ref{definition_transitive_adjacency}) to be pairwise adjacent. Similarly, the
right-neighboring links (such as $i_{-1},\cdots,i_{-4}$ in the example of
Figure \ref{link_adjacency_init} (a)) are also pairwise adjacent. This
explains the two all-$1$-triangles (surrounded by the dashed lines) in Figure
\ref{link_adjacency_init} (b), the upper one corresponding to $i^{+}$ and the
second triangle corresponding to pairwise adjacent links $i_{-1},\cdots
,i_{-4}$.

\item Since there is \textbf{at most one} link between two nodes (multi-links
are forbidden), each of the left-neighboring links can be adjacent to
\textbf{at most one} right neighboring link and vice versa. Hence in Figure
\ref{link_adjacency_init} (b), there exists \textbf{at most one} $1$-entry in
each row/column of the submatrix in yellow.
\end{enumerate}

We summarize this observation:

\begin{criterion}
If the given link adjacency pattern has the following features, it is the link
adjacency pattern of a link $i$ and its neighboring links (the labels of the
left-neighboring links are larger than link $i$, and smaller than the
right-neighboring links),\label{criterion_link_neighboring}

\begin{itemize}
\item All entries of the first row are $1$-entries;

\item The triangle bounded by the $\left(  n_{i^{+}}+1\right)  $th column
(including the $\left(  n_{i^{+}}+1\right)  $th column) is an all-$1$%
-triangle, where $n_{i^{+}}$ denotes the number of the left-neighboring links
of link $i$ and $n_{i^{+}}\geq3$;

\item There is at most one $1$-entry in each row/column of the submatrix,
which is from the $2$nd to the $\left(  n_{i^{+}}+1\right)  $th row and from
the $\left(  n_{i^{+}}+2\right)  $th to the $\left(  n_{i^{+}}+n_{i^{-}%
}+2\right)  $th column, where $n_{i^{-}}$ denotes the number of the
right-neighboring links;

\item The triangle bounded by the $\left(  n_{i^{+}}+2\right)  $th row
(including the $\left(  n_{i^{+}}+2\right)  $th row) is an all-$1$-triangle.
\end{itemize}
\end{criterion}

%

\begin{theorem}
Consider three links $i$, $j$ and $k$ are pairwise adjacent. If each of the
other $m$ links is adjacent to all the three links $i$, $j$ and $k$, then all
the $m+3$ links are pairwise adjacent. \label{theorem_relabeling_three_lines}
\end{theorem}

\begin{proof}
The three links $i$, $j$ and $k$ are pairwise adjacent and the configuration
of $i$, $j$ and $k$ can be $K_{3}$ or $K_{1,3}$, as shown in Figure
\ref{s1smaller4_2} (b). If the configuration is $K_{3}$, other links can be
adjacent to at most two of $i$, $j$ and $k$. However, if the other $m$ links
are adjacent to $i$, $j$ and $k$, the configuration of $i$, $j$ and $k$ must
be $K_{1,3}$, and $i$, $j$ and $k$ have a common endnode. Since each of the
$m$ links is adjacent to $i$, $j$ and $k$, the common endnode of $i$, $j$ and
$k$ must be also an endnode of each of the $m$ links. According to Definition
\ref{definition_transitive_adjacency}, all these $m+3$\textit{ }links are
pairwise adjacent.
\end{proof}

In Figure \ref{link_adjacency_init} (b), links $i$, $i_{+1}$ and $i_{+2}$ are
pairwise adjacent, as shown by entries in green. Links $i_{+3}$, $i_{+4}$ and
$i_{+5}$\ are adjacent to $i$, $i_{+1}$ and $i_{+2}$, as shown by entries in
magenta. By Theorem \ref{theorem_relabeling_three_lines}, links $i$, $i_{+1}$,
$i_{+2}$, $i_{+3}$, $i_{+4}$ and $i_{+5}$ are pairwise adjacent.

\subsection{The basic forbidden link adjacency patterns in a LAM}

\begin{figure}
[ptb]
\begin{center}
\includegraphics[
height=1.2713in,
width=2.3324in
]%
{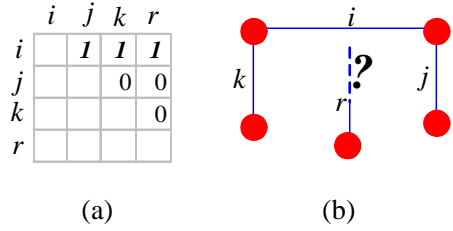}%
\caption[The smallest forbidden link adjacency pattern.]{The smallest
forbidden link adjacency pattern.}%
\label{forbidden_linkadjacency}%
\end{center}
\end{figure}

Figure \ref{forbidden_linkadjacency} (a) depicts the smallest forbidden link
adjacency pattern in a LAM. The configuration of links $i$, $j$ and $k$ is a
path on four nodes. Since link $i$ has neighboring links at both of its two
endnodes, and if link $r$ is adjacent with link $i$, then link $r$ must be
also adjacent with link $j$ or $k$. Hence, the pattern in Figure
\ref{forbidden_linkadjacency} (a)\ will not appear in a LAM.

\begin{figure}
[ptb]
\begin{center}
\includegraphics[
height=2.1897in,
width=4.1701in
]%
{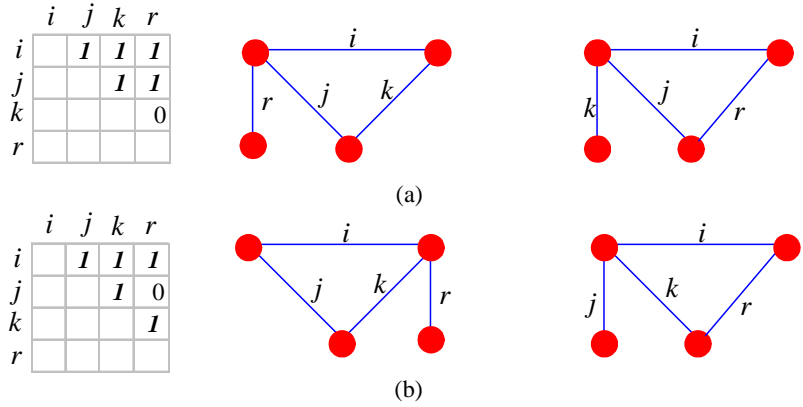}%
\caption{The possible configurations for two link adjacency patterns of $4$
links. This figure helps to prove that the patterns of $5$ links in Figure
\ref{more_forbidden_patterns} are forbidden.}%
\label{configuration_5links_aid}%
\end{center}
\end{figure}

\begin{figure}
[ptb]
\begin{center}
\includegraphics[
height=2.6091in,
width=3.4653in
]%
{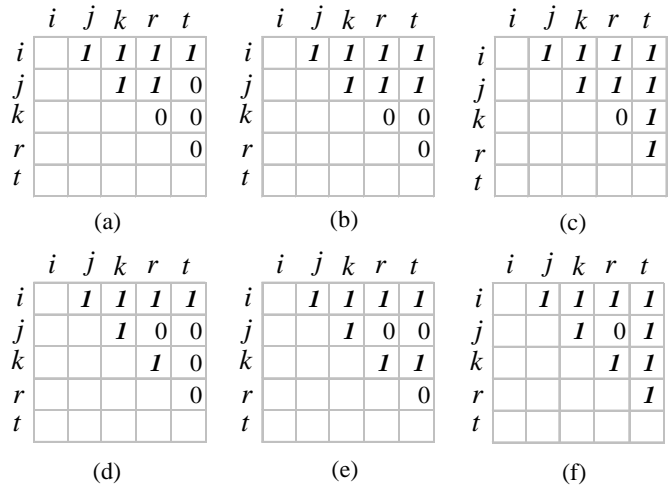}%
\caption{The forbidden link adjacency patterns of $5$ links.}%
\label{more_forbidden_patterns}%
\end{center}
\end{figure}

There are $6$ forbidden link adjacency patterns of links $i$, $j$, $k$, $r$
and $t$, as shown in Figure \ref{more_forbidden_patterns}. Since the number of
the left-neighboring links of link $i$ is smaller than $3$, we cannot use
Criterion \ref{criterion_link_neighboring} to prove that the $6$ link
adjacency patterns are forbidden. However, Figure
\ref{configuration_5links_aid}, which exhibits the possible configurations of
the link adjacency patterns of links $i$, $j$, $k$ and $r$, will facilitate
the proof that the $6$ link adjacency patterns in Figure
\ref{more_forbidden_patterns} are forbidden.

The link adjacency pattern of links $i$, $j$, $k$ and $r$ in Figure
\ref{more_forbidden_patterns} (a), (b) and (c) are the same as the link
adjacency pattern of links $i$, $j$, $k$ and $r$ in Figure
\ref{configuration_5links_aid} (a). There are only two possible configurations
of this link adjacency pattern. As we can observe in Figure
\ref{configuration_5links_aid} (a), it is impossible to have a new link $t$
which is only adjacent with link $i$, or only adjacent with links $i$ and $j$,
or adjacent with all of $i$, $j$, $k$\ and $r$. Hence, the patterns in Figure
\ref{more_forbidden_patterns} (a), (b) and (c) are forbidden. In the same way,
we observe that the patterns in Figure \ref{more_forbidden_patterns} (d), (e)
and (f) are also forbidden.%


When the number of the left-neighboring links of link $i$ is not smaller than
$3$ (which implies that the number of $1$-entries in the first all-$1$%
-triangle is not smaller than $6$), we can use Criterion
\ref{criterion_link_neighboring} to determine whether a link adjacency pattern
is forbidden.

\section{The matrix relabeling inverse line graph algorithm (MARINLINGA)}

MARINLINGA is the algorithm that we designed to compute the original graph of
a line graph, given the adjacency matrix of that line graph\footnote{Although
MARINLINGA is designed for connected line graphs, it is also convenient to
compute the original graph of a disconnected line graph component by
component. In the description of MARINLINGA, the connectedness of the
concerned graph is always assumed.}.

As explained in Section \ref{LAM_and_line_graph}, the adjacency matrix
$A_{l\left(  G\right)  }$ of $l\left(  G\right)  $ is equal to the LAM $C_{G}$
of $G$. Constructing the original graph of a line graph, is equivalent to
constructing a graph given the LAM of that graph. MARINLINGA only deals with
the upper triangle of the given LAM $C$.

\subsection{Matrix relabeling}

\label{relabeling_strategy}

\begin{figure}
[ptb]
\begin{center}
\includegraphics[
height=3.1531in,
width=6.6841in
]%
{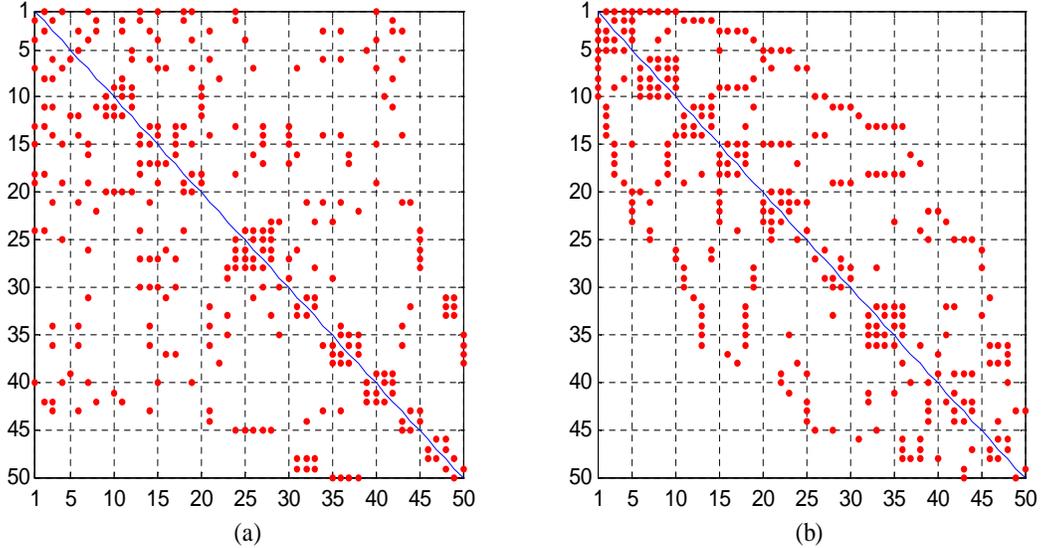}%
\caption{Matrix relabelling on the LAM $C$ of a graph with $50$ links. The red
dots represent $1$-entries. (a) Before relabelling; (b) After relabelling.}%
\label{before&after_relabelling}%
\end{center}
\end{figure}

The matrix relabeling algorithm rearranges the LAM $C$ in such a way that the
left and right neighboring links of the first link can be recognized via
Theorem \ref{theorem_relabeling_three_lines} and the construction algorithm
can work efficiently. In each column there are some $1$-entries (red dots). If
after relabeling the top $1$-entries of all the columns are connected by a
curve, the curve should be nonincreasing. For example, by the LAM $C$ of a
graph with $50$ links in Figure \ref{before&after_relabelling} (a), we can
only determine which links are adjacent to the first link, without any
information about which endnode of the first link that the neighboring links
are incident to. Fortunately, according to Theorem
\ref{theorem_relabeling_three_lines}, the relabeled LAM $C$ in Figure
\ref{before&after_relabelling} (b) tells that links $2$-$5$ are the
left-neighboring links of the first link and links $6$-$10$ are the
right-neighboring links.%

\begin{figure}
[ptb]
\begin{center}
\includegraphics[
height=2.3791in,
width=4.2722in
]%
{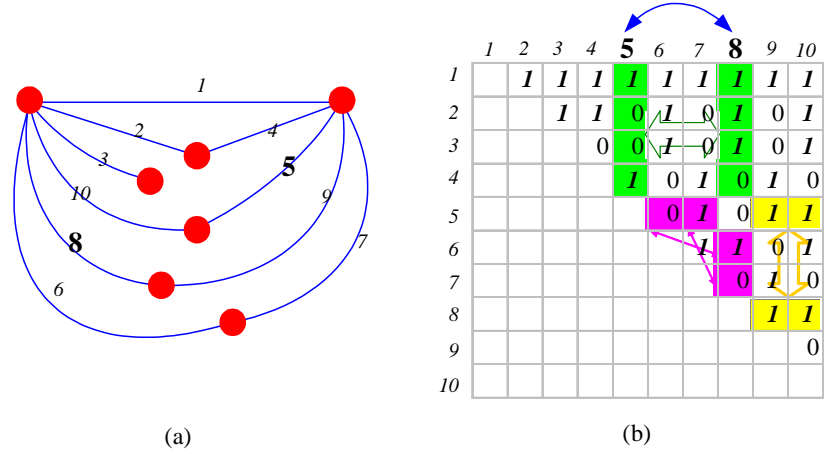}%
\caption{The illustration of swapping the labels of link $5$ and $8$. The
entries in green, magenta and yellow ought to be swapped respectively.}%
\label{label_swapping}%
\end{center}
\end{figure}

Let us first introduce the meaning of swapping the labels of two links in a
LAM $C_{L_{G}\times L_{G}}$. The entry $c_{ij}$ indicates whether links $i$
and $j$ are adjacent. Swapping the labels of links $j$ and $k$ ($j<k$) implies
that links which are previously adjacent to link $j$ are now adjacent to link
$k$, and links which are previously adjacent to link $k$, are now adjacent to
link $j$, but the adjacency relation between links $j$ and $k$ is the same as
before, namely the entry $c_{jk}$ of $C_{L_{G}\times L_{G}}$ is unchanged.
Hence, swapping the labels of links $j$ and $k$ ($j<k)$ means to swap the
entries $c_{ij}$ and $c_{ik}$ for $i=1,2,\cdots,j-1$ (shown in the example of
Figure \ref{label_swapping} in green), the entries $c_{ji}$ and $c_{ik}$ for
$i=j+1,\cdots,k-1$ (in magenta), the entries $c_{ji}$ and $c_{ki}$,
$i=k+1,\cdots,L_{G}-1,L_{G}$ (in yellow).%

\begin{algorithm}%

\caption{$C \Leftarrow SwapLabel(C,j,k)$}\label{algo_label_swapping}

\begin{algorithmic}[1]
\FOR{$i=1$ to $j-1$}
\STATE $swap(c_{ij},c_{ik})$
\ENDFOR
\FOR{$i=j+1$ to $k-1$}
\STATE $swap(c_{ji},c_{ik})$
\ENDFOR
\FOR{$i=k+1$ to $L_{G}$}
\STATE $swap(c_{ji},c_{ki})$
\ENDFOR
\end{algorithmic}%

\end{algorithm}%

Lines 1-3 of the metacode of Algorithm \ref{algo_label_swapping} swap the
entries $c_{ji}$ and $c_{ik}$, $i=j+1,\cdots,k-1$, and lines 4-6 swap the
entries $c_{ji}$ and $c_{ik}$, $i=j+1,\cdots,k-1$, and lines 7-9 swap the
entries $c_{ji}$ and $c_{ki}$, $i=k+1,\cdots,L_{G}-1,L_{G}$. The code
$swap\left(  c_{ij},c_{ik}\right)  $ of line 2 is equivalent to the codes:
$t=c_{ij}$; $c_{ij}=c_{ik}$; $c_{ik}=t$.%

\begin{figure}
[ptb]
\begin{center}
\includegraphics[
height=3.5241in,
width=4.6855in
]%
{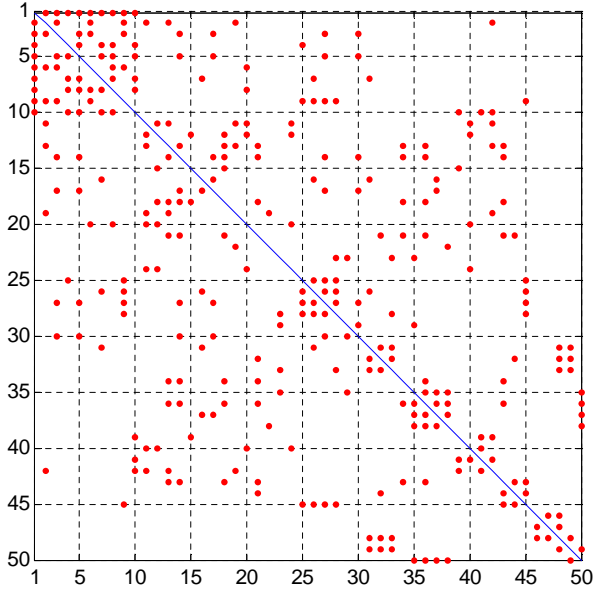}%
\caption{The LAM $C$ after the relabeling of the first row.}%
\label{example_relabeling}%
\end{center}
\end{figure}

Next, we will explain the matrix relabeling algorithm. We will first give an
example showing how the matrix relabeling algorithm relabels the LAM $C$ in
Figure \ref{before&after_relabelling} (a) into the matrix in Figure
\ref{before&after_relabelling} (b). In the first row of the matrix in Figure
\ref{before&after_relabelling} (a) there are $9$ 1-entries in total. There are
$6$ 0-entries from $c_{1,2}$ to $c_{1,10}$ and $6$ 1-entries from $c_{1,11}$
to $c_{1,50}$: $c_{1,3}=c_{1,5}=c_{1,6}=c_{1,8}=c_{1,9}=c_{1,10}=0$ and
$c_{1,13}=c_{1,15}=c_{1,18}=c_{1,19}=c_{1,24}=c_{1,40}=1$. We swap the labels
of links $3$ and $13$, links $5$ and $15$, links $6$ and $18$, links $8$ and
$19$, links $9$ and $24$, links $10$ and $40$ by Algorithm
\ref{algo_group_label_swapping} and the LAM $C$ is shown in Figure
\ref{example_relabeling}. In the second row, there are $3$ 1-entries from
$c_{2,3}$ to $c_{2,10}$. There are $2$ 0-entries from $c_{2,3}$ to $c_{2,5}$
and $2$ 1-entries from $c_{2,6}$ to $c_{2,10}$: $c_{2,4}=c_{2,5}=0$ and
$c_{2,6}=c_{2,9}=1$. We swap the labels of links $4$ and $6$, links $5$ and
$9$. By similar operations, we relabel the LAM $C$ into the order shown in
Figure \ref{before&after_relabelling} (b).
%

\begin{algorithm}%

\caption{$C \Leftarrow GroupLabelSwapping(C,u,k,a,b)$}
\label{algo_group_label_swapping}

\begin{algorithmic}[1]
\STATE $m \Leftarrow 0$
\FOR{$i=u+1$ to $a+u$}
\IF{$c_{ki}=0$}
\STATE $m \Leftarrow m+1$
\STATE $X_{m} \Leftarrow i$
\ENDIF
\ENDFOR
\STATE $m \Leftarrow$ $0$
\FOR{$i=a+u+1$ to $b$}
\IF{$c_{ki}=1$}
\STATE $m \Leftarrow m+1$
\STATE $Y_{m} \Leftarrow i$
\ENDIF
\ENDFOR
\FOR{$i=1$ to $m$}
\STATE $C \Leftarrow SwapLabel(C,X_{i},Y_{i})$
\ENDFOR
\end{algorithmic}%

\end{algorithm}%

Now we give the general description of the matrix relabeling algorithm. In the
$k$th row of $C$, Lines 1-7 of Algorithm \ref{algo_group_label_swapping} store
the value of $i$ in $X$ when the entry $c_{ki}$ is $0$, $i=u+1,\cdots,a+u$.
Lines 8-14 store the value of $i$ in $Y$ when the entry $c_{ki}$ is $1$,
$i=a+u+1,\cdots,b$. If $a=%
{\textstyle\sum\limits_{i=u+1}^{b}}
c_{ki}$, $X$ and $Y$ have the same number of elements. Lines 15-17 swap the
labels of $X_{i}$ and $Y_{i}$, where $X_{i}$ and $Y_{i}$ are the $i$th element
of $X$ and $Y$ respectively. For instance in Figure \ref{label_swapping} (b),
if we take $u=2$, $k=2$, $b=10$ and $a=%
{\textstyle\sum\limits_{i=3}^{10}}
c_{2i}$ $=5$, by Algorithm \ref{algo_group_label_swapping}, $X=\left[
\begin{array}
[c]{cc}%
5 & 7
\end{array}
\right]  ^{T},Y=\left[
\begin{array}
[c]{cc}%
8 & 10
\end{array}
\right]  ^{T}$, the labels of links $5$ and $8$, $7$ and $10$ are swapped respectively.%

\begin{algorithm}%

\caption{$(C,s_{1},s_{2},s_{3}) \Leftarrow MatrixRelabeling(C)$}
\label{algo_link_relabeling}

\begin{algorithmic}[1]
\STATE $s_{1} \Leftarrow$ the sum of $c_{1i}$, where $i=2$ to $L_{G}$
\STATE $C \Leftarrow GroupLabelSwapping(C,1,1,s_{1},L_{G})$
\STATE $s_{2} \Leftarrow$ the sum of $c_{2i}$, where $i=3$ to $s_{1}+1$
\STATE $C \Leftarrow GroupLabelSwapping(C,2,2,s_{2},s_{1}+1)$
\STATE $s_{3} \Leftarrow$ the sum of $c_{3i}$, where $i=4$ to $s_{2}+2$
\STATE $C \Leftarrow GroupLabelSwapping(C,3,3,s_{3},s_{2}+2)$
\STATE $ \bar{s} \Leftarrow s_{1}+1$
\STATE $ k \Leftarrow 2$
\WHILE{$ \bar{s} < L_{G}$ and $k \leq L_{G}$}
\STATE $s \Leftarrow$ the sum of $c_{ki}$, where $i=\bar{s}+1$ to $L_{G}$
\STATE $C \Leftarrow GroupLabelSwapping(C,\bar{s},k,s,L_{G})$
\STATE $ k \Leftarrow k+1$
\STATE $ \bar{s} \Leftarrow \bar{s}+s$
\ENDWHILE
\end{algorithmic}%

\end{algorithm}%

\begin{figure}
[ptb]
\begin{center}
\includegraphics[
height=5.0488in,
width=5.1387in
]%
{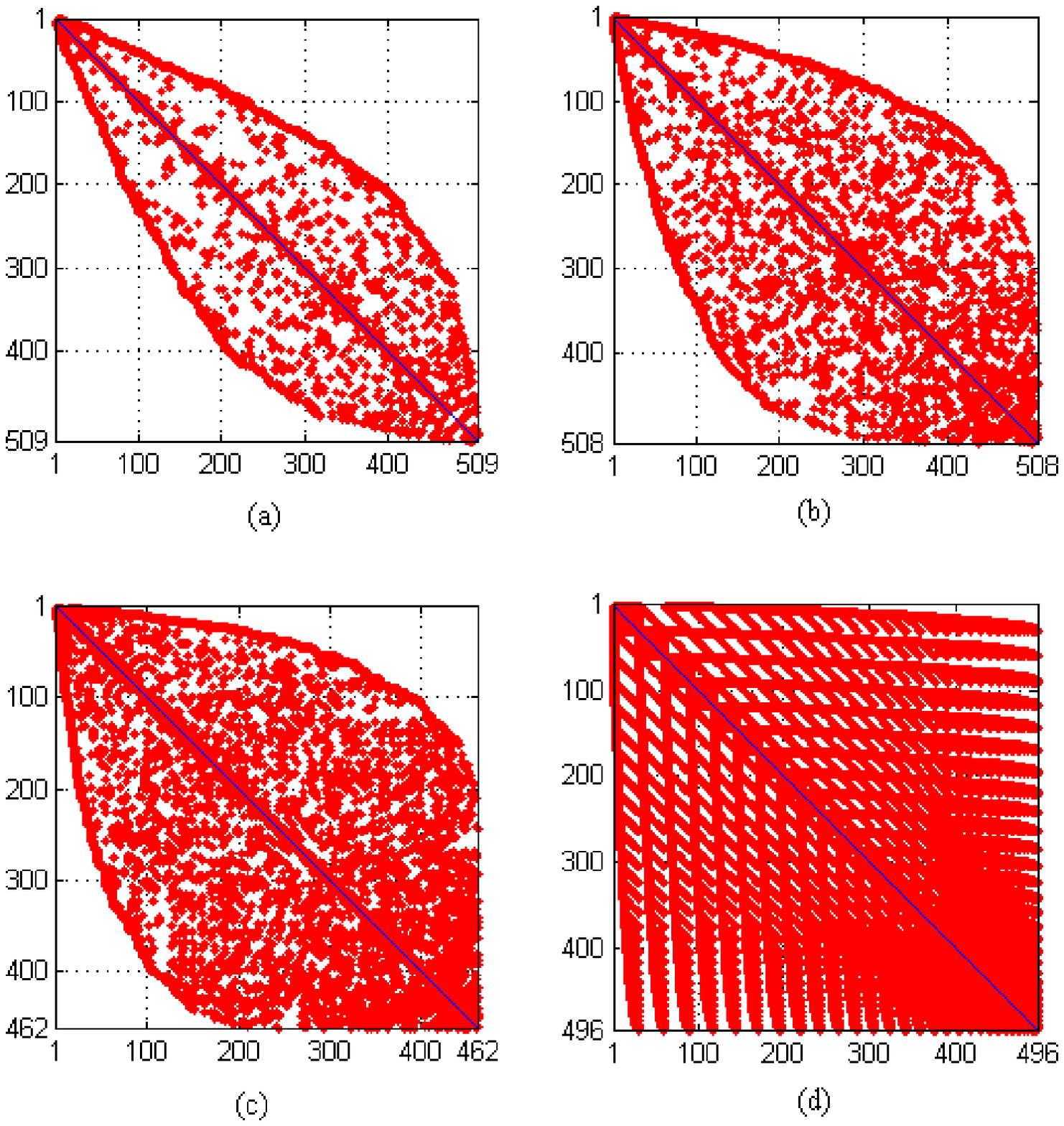}%
\caption{The relabeled $C$ of four ER random graphs $G\left(  N,p\right)  $:
(a) $N=350$, $p=\frac{\log\left(  N\right)  }{2N}$; (b) $N=200$, $p=\frac
{\log\left(  N\right)  }{N}$; (c) $N=100$, $p=\frac{2\log\left(  N\right)
}{N}$; (d) $N=32$, $p=1$, where $p=\frac{\log\left(  N\right)  }{N}$ is the
threshold probability for the connectivity of the graph.}%
\label{examples_relabeling}%
\end{center}
\end{figure}

Lines 1-2 of Algorithm \ref{algo_link_relabeling} make the neighboring links
of link $1$ have the smaller labels than the other links. By lines 3-4, the
labels of the links which are adjacent to both link $1$ and $2$ are smaller
than those of the remaining links. Further, lines 5-6 let the labels of the
links which are adjacent to all of links $1$, $2$ and $3$ are smaller than
those of the remaining links. Lines 7-14 make that the labels of the links
which are adjacent to link $i$ but not adjacent to links $1,\cdots,i-1$, are
smaller than the labels of the links which are not adjacent to link
$1,\cdots,i$, for $i=2,\cdots,L_{G}$. Figure \ref{before&after_relabelling}
and \ref{examples_relabeling} show examples of $C$ before and after matrix
relabeling.%

Let $s_{1}=%
{\textstyle\sum\limits_{i=2}^{L_{G}}}
c_{1i}$, $s_{2}=%
{\textstyle\sum\limits_{i=3}^{s_{1}+1}}
c_{2i}$ and $s_{3}=%
{\textstyle\sum\limits_{i=4}^{s_{2}+2}}
c_{3i}$. After relabeling by Algorithm \ref{algo_link_relabeling},\ the given
LAM $C$ satisfies:

\begin{itemize}
\item For $i=2,\cdots,s_{1}+1$, $c_{1i}=1$; and for $i=s_{1}+2,\cdots,L_{G}$,
$c_{1i}=0$.

\item For $i=3,\cdots,s_{2}+2$, $c_{2i}=1$ if $s_{2}\geq1$; and for
$i=s_{2}+3,\cdots,s_{1}+1$, $c_{2i}=0$ if $s_{1}\geq s_{2}+2$.

\item For $i=4,\cdots,s_{3}+3$, $c_{3i}=1$ if $s_{3}\geq1$; and for
$i=s_{3}+4,\cdots,s_{2}+2$, $c_{3i}=0$ if $s_{2}\geq s_{3}+2$.

\item If link $j$ ($j\geq s_{1+1}$) is adjacent to link $i$ but not adjacent
to links $1,2,\cdots,i-1$ ($i\geq2$), and link $k$ ($k\geq s_{1+1}$) is not
adjacent to all of links $1,2,\cdots,i$ ($i\geq2$), then $j<k$.
\end{itemize}

\begin{figure}
[ptb]
\begin{center}
\includegraphics[
height=2.2217in,
width=4.2575in
]%
{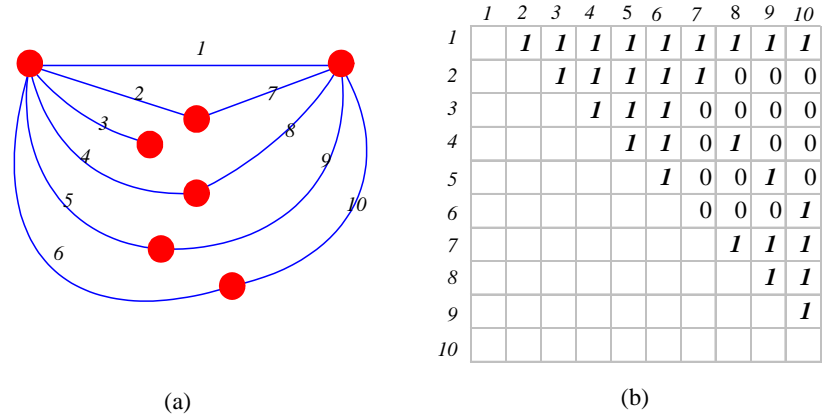}%
\caption{The LAM (a) relabeled by Algorithm \ref{algo_link_relabeling} and its
corresponding graph (b).}%
\label{relabeled_example}%
\end{center}
\end{figure}

If $s_{3}\geq1$ (which implies that $s_{2}\geq2$ and $s_{1}\geq3$), according
to Theorem \ref{theorem_relabeling_three_lines}, links $2,3,\cdots,s_{3}+3$
are the left-neighboring links of links $1$ and the links $s_{3}%
+4,\cdots,s_{1}+1$ are the right-neighboring links of link $1$, as illustrated
in the example of Figure \ref{relabeled_example} where $s_{1}=9$ and $s_{3}%
=3$.%

\subsection{Construction algorithm}

The construction algorithm converts the relabeled $C$ into the matrix
$E_{2\times L_{G}}$, where the entries $e_{1i}$ and $e_{2i}$ denotes the two
endnodes of link $i$. During the process of the construction, the zero entries
of $E_{2\times L_{G}}$ mean that the endnodes have not been determined yet.

Section \ref{Section_construction_algo_example} will first show an example of
graph construction, and section \ref{Section_construction_algo_intial} and
\ref{Section_construction_algo_} will describe the general construction algorithm.

\subsubsection{An example of graph construction from $C$}

\label{Section_construction_algo_example}From the given LAM $C$ in Figure
\ref{before&after_relabelling} (b), we deduce that the graph has $50$ links.
Based on the LAM $C$, we will determine the endnodes of the $50$ links. The
construction consists of the following steps:\noindent

\begin{enumerate}
\item Let nodes $1$ and $2$ be the endnodes of link $1$. According to Theorem
\ref{theorem_relabeling_three_lines}, node $1$ is also the endnode of links
$2$-$5$ and node $2$ is also the endnode of links $6$-$10$, as shown in Figure
\ref{1st_phase} (a) and equation $\left(  \ref{Eq_construction_example_init}%
\right)  $ below, where the numbers above the matrix are the link numbers.%
\begin{equation}%
\begin{array}
[c]{c}%
\\
E=
\end{array}%
\begin{array}
[c]{c}%
\\
\left[
\begin{array}
[c]{c}%
\\
\\
\end{array}
\right.
\end{array}%
\begin{tabular}
[c]{ccccccccccccc}%
${\footnotesize 1}$ & ${\footnotesize 2}$ & ${\footnotesize 3}$ &
${\footnotesize 4}$ & ${\footnotesize 5}$ & ${\footnotesize 6}$ &
${\footnotesize 7}$ & ${\footnotesize 8}$ & ${\footnotesize 9}$ &
${\footnotesize 10}$ & ${\footnotesize 11}$ & ${\footnotesize \cdots}$ &
${\footnotesize 50}$\\\hline
\multicolumn{1}{|c}{$1$} & \multicolumn{1}{|c}{$1$} & \multicolumn{1}{|c}{$1$}
& \multicolumn{1}{|c}{$1$} & \multicolumn{1}{|c}{$1$} &
\multicolumn{1}{|c}{$2$} & \multicolumn{1}{|c}{$2$} & \multicolumn{1}{|c}{$2$}
& \multicolumn{1}{|c}{$2$} & \multicolumn{1}{|c}{$2$} &
\multicolumn{1}{|c}{$0$} & \multicolumn{1}{|c}{$\cdots$} &
\multicolumn{1}{|c|}{$0$}\\\hline
\multicolumn{1}{|c}{$2$} & \multicolumn{1}{|c}{$0$} & \multicolumn{1}{|c}{$0$}
& \multicolumn{1}{|c}{$0$} & \multicolumn{1}{|c}{$0$} &
\multicolumn{1}{|c}{$0$} & \multicolumn{1}{|c}{$0$} & \multicolumn{1}{|c}{$0$}
& \multicolumn{1}{|c}{$0$} & $0$ & \multicolumn{1}{|c}{$0$} &
\multicolumn{1}{|c}{$\cdots$} & \multicolumn{1}{|c|}{$0$}\\\hline
\end{tabular}%
\begin{array}
[c]{c}%
\\
\left.
\begin{array}
[c]{c}%
\\
\\
\end{array}
\right]
\end{array}
\label{Eq_construction_example_init}%
\end{equation}
Let node $3$ be the other endnode of link $2$. The $2$nd row of the LAM $C$
shows that links $11$-$14$ are adjacent to link $2$. Hence, node $3$ is also
the endnode of links $11$-$14$, as shown in Figure \ref{1st_phase} (b) and
equation $\left(  \ref{Eq_construction_example_1_1}\right)  $.%
\begin{equation}%
\begin{array}
[c]{c}%
\\
E=
\end{array}%
\begin{array}
[c]{c}%
\\
\left[
\begin{array}
[c]{c}%
\\
\\
\end{array}
\right.
\end{array}%
\begin{tabular}
[c]{ccccccccccccccccc}%
${\footnotesize 1}$ & ${\footnotesize 2}$ & ${\footnotesize 3}$ &
${\footnotesize 4}$ & ${\footnotesize 5}$ & ${\footnotesize 6}$ &
${\footnotesize 7}$ & ${\footnotesize 8}$ & ${\footnotesize 9}$ &
${\footnotesize 10}$ & ${\footnotesize 11}$ & ${\footnotesize 12}$ &
${\footnotesize 13}$ & ${\footnotesize 14}$ & ${\footnotesize 15}$ &
${\footnotesize \cdots}$ & ${\footnotesize 50}$\\\hline
\multicolumn{1}{|c}{$1$} & \multicolumn{1}{|c}{$1$} & \multicolumn{1}{|c}{$1$}
& \multicolumn{1}{|c}{$1$} & \multicolumn{1}{|c}{$1$} &
\multicolumn{1}{|c}{$2$} & \multicolumn{1}{|c}{$2$} & \multicolumn{1}{|c}{$2$}
& \multicolumn{1}{|c}{$2$} & \multicolumn{1}{|c}{$2$} &
\multicolumn{1}{|c}{$3$} & \multicolumn{1}{|c}{$3$} & \multicolumn{1}{|c}{$3$}
& \multicolumn{1}{|c}{$3$} & \multicolumn{1}{|c}{$0$} &
\multicolumn{1}{|c}{$\cdots$} & \multicolumn{1}{|c|}{$0$}\\\hline
\multicolumn{1}{|c}{$2$} & \multicolumn{1}{|c}{$3$} & \multicolumn{1}{|c}{$0$}
& \multicolumn{1}{|c}{$0$} & \multicolumn{1}{|c}{$0$} &
\multicolumn{1}{|c}{$0$} & \multicolumn{1}{|c}{$0$} & \multicolumn{1}{|c}{$0$}
& \multicolumn{1}{|c}{$0$} & \multicolumn{1}{|c}{$0$} &
\multicolumn{1}{|c}{$0$} & \multicolumn{1}{|c}{$0$} & \multicolumn{1}{|c}{$0$}
& \multicolumn{1}{|c}{$0$} & \multicolumn{1}{|c}{$0$} &
\multicolumn{1}{|c}{$\cdots$} & \multicolumn{1}{|c|}{$0$}\\\hline
\end{tabular}%
\begin{array}
[c]{c}%
\\
\left.
\begin{array}
[c]{c}%
\\
\\
\end{array}
\right]
\end{array}
\label{Eq_construction_example_1_1}%
\end{equation}
Similarly, let node $4$ be the endnode of link $3$, $6$ and $15$-$18$ as shown
in Figure \ref{1st_phase} (c) and equation $\left(
\ref{Eq_construction_example_1_2}\right)  $,%
\begin{equation}%
\begin{array}
[c]{c}%
\\
E=
\end{array}%
\begin{array}
[c]{c}%
\\
\left[
\begin{array}
[c]{c}%
\\
\\
\end{array}
\right.  \\
\\
\\
\\
\end{array}%
\begin{tabular}
[c]{cccccccccccccc}%
${\footnotesize 1}$ & ${\footnotesize 2}$ & ${\footnotesize 3}$ &
${\footnotesize 4}$ & ${\footnotesize 5}$ & ${\footnotesize 6}$ &
${\footnotesize 7}$ & ${\footnotesize 8}$ & ${\footnotesize 9}$ &
${\footnotesize 10}$ & ${\footnotesize 11}$ & ${\footnotesize 12}$ &
${\footnotesize 13}$ & ${\footnotesize 14}$\\\hline
\multicolumn{1}{|c}{$1$} & \multicolumn{1}{|c}{$1$} & \multicolumn{1}{|c}{$1$}
& \multicolumn{1}{|c}{$1$} & \multicolumn{1}{|c}{$1$} &
\multicolumn{1}{|c}{$2$} & \multicolumn{1}{|c}{$2$} & \multicolumn{1}{|c}{$2$}
& \multicolumn{1}{|c}{$2$} & \multicolumn{1}{|c}{$2$} &
\multicolumn{1}{|c}{$3$} & \multicolumn{1}{|c}{$3$} & \multicolumn{1}{|c}{$3$}
& \multicolumn{1}{|c|}{$3$}\\\hline
\multicolumn{1}{|c}{$2$} & \multicolumn{1}{|c}{$3$} & \multicolumn{1}{|c}{$4$}
& \multicolumn{1}{|c}{$0$} & \multicolumn{1}{|c}{$0$} &
\multicolumn{1}{|c}{$4$} & \multicolumn{1}{|c}{$0$} & \multicolumn{1}{|c}{$0$}
& \multicolumn{1}{|c}{$0$} & \multicolumn{1}{|c}{$0$} &
\multicolumn{1}{|c}{$0$} & \multicolumn{1}{|c}{$0$} & \multicolumn{1}{|c}{$0$}
& \multicolumn{1}{|c|}{$0$}\\\hline
${\footnotesize 15}$ & ${\footnotesize 16}$ & ${\footnotesize 17}$ &
${\footnotesize 18}$ & ${\footnotesize 19}$ & ${\footnotesize \cdots}$ &
${\footnotesize 50}$ &  &  &  &  &  &  & \\\cline{1-7}%
\multicolumn{1}{|c}{$4$} & \multicolumn{1}{|c}{$4$} & \multicolumn{1}{|c}{$4$}
& \multicolumn{1}{|c}{$4$} & \multicolumn{1}{|c}{$0$} &
\multicolumn{1}{|c}{$\cdots$} & \multicolumn{1}{|c}{$0$} &
\multicolumn{1}{|c}{} &  &  &  &  &  & \\\cline{1-7}%
\multicolumn{1}{|c}{$0$} & \multicolumn{1}{|c}{$0$} & \multicolumn{1}{|c}{$0$}
& \multicolumn{1}{|c}{$0$} & \multicolumn{1}{|c}{$0$} &
\multicolumn{1}{|c}{$\cdots$} & \multicolumn{1}{|c}{$0$} &
\multicolumn{1}{|c}{} &  &  &  &  &  & \\\cline{1-7}%
\end{tabular}%
\begin{array}
[c]{c}%
\\
\\
\\
\\
\left.
\begin{array}
[c]{c}%
\\
\\
\end{array}
\right]  \\
\end{array}
\label{Eq_construction_example_1_2}%
\end{equation}
and let node $5$ be the endnode of link $4$, $8$ and $19$ as shown in Figure
\ref{1st_phase} (d) and equation $\left(  \ref{Eq_construction_example_1_3}%
\right)  $,%
\begin{equation}%
\begin{array}
[c]{c}%
\\
E=
\end{array}%
\begin{array}
[c]{c}%
\\
\left[
\begin{array}
[c]{c}%
\\
\\
\end{array}
\right.  \\
\\
\\
\\
\end{array}%
\begin{tabular}
[c]{cccccccccccccc}%
${\footnotesize 1}$ & ${\footnotesize 2}$ & ${\footnotesize 3}$ &
${\footnotesize 4}$ & ${\footnotesize 5}$ & ${\footnotesize 6}$ &
${\footnotesize 7}$ & ${\footnotesize 8}$ & ${\footnotesize 9}$ &
${\footnotesize 10}$ & ${\footnotesize 11}$ & ${\footnotesize 12}$ &
${\footnotesize 13}$ & ${\footnotesize 14}$\\\hline
\multicolumn{1}{|c}{$1$} & \multicolumn{1}{|c}{$1$} & \multicolumn{1}{|c}{$1$}
& \multicolumn{1}{|c}{$1$} & \multicolumn{1}{|c}{$1$} &
\multicolumn{1}{|c}{$2$} & \multicolumn{1}{|c}{$2$} & \multicolumn{1}{|c}{$2$}
& \multicolumn{1}{|c}{$2$} & \multicolumn{1}{|c}{$2$} &
\multicolumn{1}{|c}{$3$} & \multicolumn{1}{|c}{$3$} & \multicolumn{1}{|c}{$3$}
& \multicolumn{1}{|c|}{$3$}\\\hline
\multicolumn{1}{|c}{$2$} & \multicolumn{1}{|c}{$3$} & \multicolumn{1}{|c}{$4$}
& \multicolumn{1}{|c}{$5$} & \multicolumn{1}{|c}{$0$} &
\multicolumn{1}{|c}{$4$} & \multicolumn{1}{|c}{$0$} & \multicolumn{1}{|c}{$5$}
& \multicolumn{1}{|c}{$0$} & \multicolumn{1}{|c}{$0$} &
\multicolumn{1}{|c}{$0$} & \multicolumn{1}{|c}{$0$} & \multicolumn{1}{|c}{$0$}
& \multicolumn{1}{|c|}{$0$}\\\hline
${\footnotesize 15}$ & ${\footnotesize 16}$ & ${\footnotesize 17}$ &
${\footnotesize 18}$ & ${\footnotesize 19}$ & ${\footnotesize 20}$ &
${\footnotesize \cdots}$ & ${\footnotesize 50}$ &  &  &  &  &  & \\\cline{1-8}%
\multicolumn{1}{|c}{$4$} & \multicolumn{1}{|c}{$4$} & \multicolumn{1}{|c}{$4$}
& \multicolumn{1}{|c}{$4$} & \multicolumn{1}{|c}{$5$} &
\multicolumn{1}{|c}{$0$} & \multicolumn{1}{|c}{$\cdots$} &
\multicolumn{1}{|c}{$0$} & \multicolumn{1}{|c}{} &  &  &  &  & \\\cline{1-8}%
\multicolumn{1}{|c}{$0$} & \multicolumn{1}{|c}{$0$} & \multicolumn{1}{|c}{$0$}
& \multicolumn{1}{|c}{$0$} & \multicolumn{1}{|c}{$0$} &
\multicolumn{1}{|c}{$0$} & \multicolumn{1}{|c}{$\cdots$} &
\multicolumn{1}{|c}{$0$} & \multicolumn{1}{|c}{} &  &  &  &  & \\\cline{1-8}%
\end{tabular}%
\begin{array}
[c]{c}%
\\
\\
\\
\\
\left.
\begin{array}
[c]{c}%
\\
\\
\end{array}
\right]  \\
\end{array}
\label{Eq_construction_example_1_3}%
\end{equation}
and let node $6$ be the endnode of link $5$, $16$ and $20$-$23$ as shown in
Figure \ref{1st_phase} (e) and equation $\left(
\ref{Eq_construction_example_1_4}\right)  $.%
\begin{equation}%
\begin{array}
[c]{c}%
\\
E=
\end{array}%
\begin{array}
[c]{c}%
\\
\left[
\begin{array}
[c]{c}%
\\
\\
\end{array}
\right.  \\
\\
\\
\\
\end{array}%
\begin{tabular}
[c]{cccccccccccccc}%
${\footnotesize 1}$ & ${\footnotesize 2}$ & ${\footnotesize 3}$ &
${\footnotesize 4}$ & ${\footnotesize 5}$ & ${\footnotesize 6}$ &
${\footnotesize 7}$ & ${\footnotesize 8}$ & ${\footnotesize 9}$ &
${\footnotesize 10}$ & ${\footnotesize 11}$ & ${\footnotesize 12}$ &
${\footnotesize 13}$ & ${\footnotesize 14}$\\\hline
\multicolumn{1}{|c}{$1$} & \multicolumn{1}{|c}{$1$} & \multicolumn{1}{|c}{$1$}
& \multicolumn{1}{|c}{$1$} & \multicolumn{1}{|c}{$1$} &
\multicolumn{1}{|c}{$2$} & \multicolumn{1}{|c}{$2$} & \multicolumn{1}{|c}{$2$}
& \multicolumn{1}{|c}{$2$} & \multicolumn{1}{|c}{$2$} &
\multicolumn{1}{|c}{$3$} & \multicolumn{1}{|c}{$3$} & \multicolumn{1}{|c}{$3$}
& \multicolumn{1}{|c|}{$3$}\\\hline
\multicolumn{1}{|c}{$2$} & \multicolumn{1}{|c}{$3$} & \multicolumn{1}{|c}{$4$}
& \multicolumn{1}{|c}{$5$} & \multicolumn{1}{|c}{$6$} &
\multicolumn{1}{|c}{$4$} & \multicolumn{1}{|c}{$0$} & \multicolumn{1}{|c}{$5$}
& \multicolumn{1}{|c}{$0$} & \multicolumn{1}{|c}{$0$} &
\multicolumn{1}{|c}{$0$} & \multicolumn{1}{|c}{$0$} & \multicolumn{1}{|c}{$0$}
& \multicolumn{1}{|c|}{$0$}\\\hline
${\footnotesize 15}$ & ${\footnotesize 16}$ & ${\footnotesize 17}$ &
${\footnotesize 18}$ & ${\footnotesize 19}$ & ${\footnotesize 20}$ &
${\footnotesize 21}$ & ${\footnotesize 22}$ & ${\footnotesize 23}$ &
${\footnotesize 24}$ & ${\footnotesize \cdots}$ & ${\footnotesize 50}$ &  &
\\\cline{1-12}\cline{7-12}%
\multicolumn{1}{|c}{$4$} & \multicolumn{1}{|c}{$4$} & \multicolumn{1}{|c}{$4$}
& \multicolumn{1}{|c}{$4$} & \multicolumn{1}{|c}{$5$} &
\multicolumn{1}{|c}{$6$} & \multicolumn{1}{|c}{$6$} & \multicolumn{1}{|c}{$6$}
& \multicolumn{1}{|c}{$6$} & \multicolumn{1}{|c}{$0$} &
\multicolumn{1}{|c}{$\cdots$} & \multicolumn{1}{|c}{$0$} &
\multicolumn{1}{|c}{} & \\\cline{1-12}\cline{7-12}%
\multicolumn{1}{|c}{$0$} & \multicolumn{1}{|c}{$6$} & \multicolumn{1}{|c}{$0$}
& \multicolumn{1}{|c}{$0$} & \multicolumn{1}{|c}{$0$} &
\multicolumn{1}{|c}{$0$} & \multicolumn{1}{|c}{$0$} & \multicolumn{1}{|c}{$0$}
& \multicolumn{1}{|c}{$0$} & \multicolumn{1}{|c}{$0$} &
\multicolumn{1}{|c}{$\cdots$} & \multicolumn{1}{|c}{$0$} &
\multicolumn{1}{|c}{} & \\\cline{1-12}\cline{7-12}%
\end{tabular}%
\begin{array}
[c]{c}%
\\
\\
\\
\\
\left.
\begin{array}
[c]{c}%
\\
\\
\end{array}
\right]  \\
\end{array}
\label{Eq_construction_example_1_4}%
\end{equation}
Then compute the LAM of the constructed part of the graph as shown in Figure
\ref{phase_1st_matrix}. The red dots are $1$-entries which are from the given
LAM\ in Figure \ref{before&after_relabelling} (b). The green dots are
$1$-entries which are determined by the red $1$-entries. If the corresponding
entries in the given matrix are not $1$, then the matrix is not a LAM.%
\begin{figure}
[ptb]
\begin{center}
\includegraphics[
height=5.1076in,
width=6.0476in
]%
{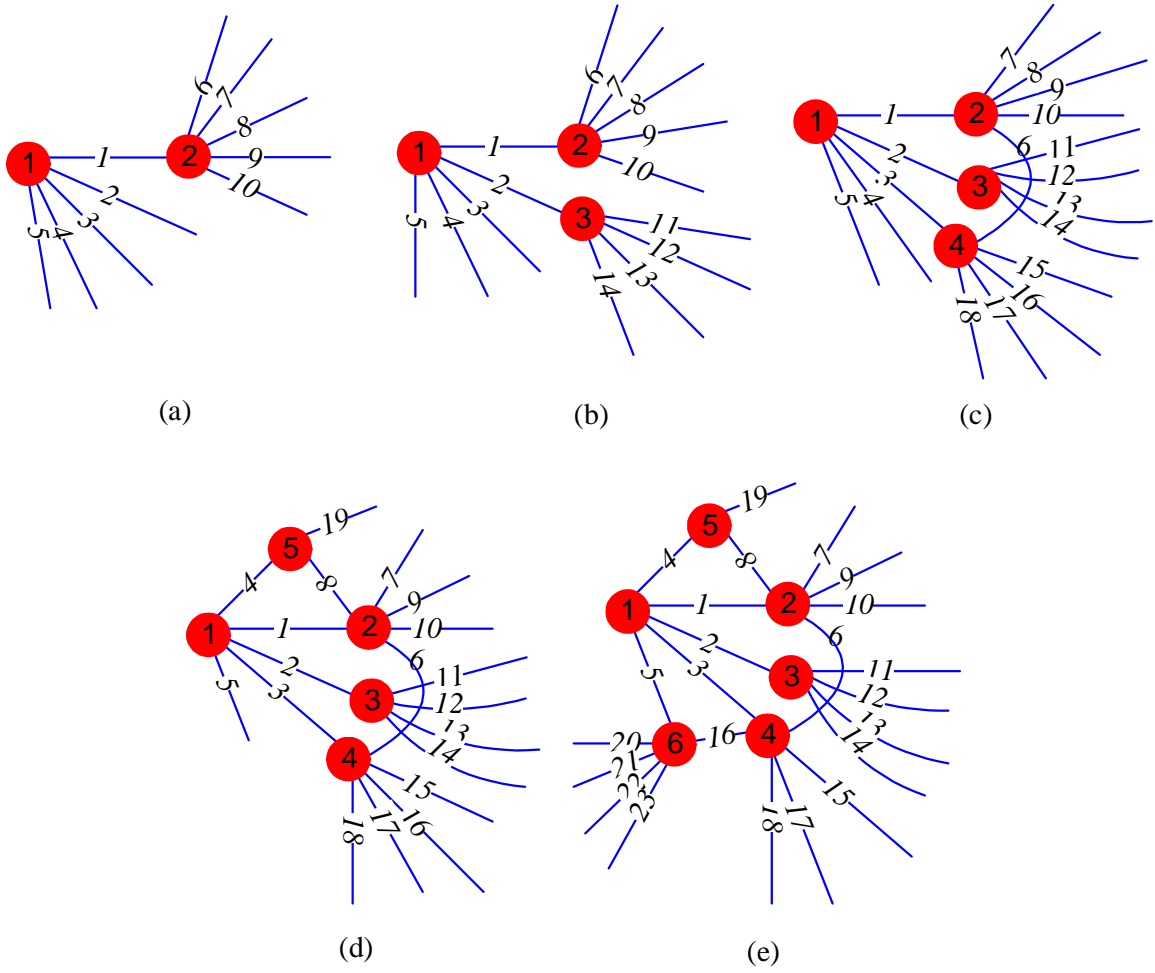}%
\caption{The example of construction. The initialization is done in (a). Both
or one of the two endnodes of links $1$-$23$ are determined.}%
\label{1st_phase}%
\end{center}
\end{figure}
\begin{figure}
[ptb]
\begin{center}
\includegraphics[
height=2.1715in,
width=2.5806in
]%
{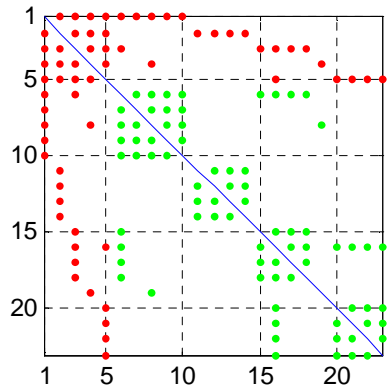}%
\caption{The LAM\ of the constructed part (links $1$-$23$) of graph are
computed. The green $1$-entries are determined by the red $1$-entries.}%
\label{phase_1st_matrix}%
\end{center}
\end{figure}

\item In the second step, we scan rows $6$ to $10$ of the LAM, since links $6$
to $10$ are incident to the same endnode. Let node $7$ be the endnode of link
$7$, $21$ and $24$-$25$, and let node $8$ be the endnode of link $9$ and $20$,
and let node $9$ be the endnode of link $10$, $14$ and $26$-$27$, as shown in
Equation $\left(  \ref{Eq_construction_example_2}\right)  $ and Figure
\ref{2nd_phase}.%
\begin{equation}%
\begin{array}
[c]{c}%
\\
E=
\end{array}%
\begin{array}
[c]{c}%
\\
\left[
\begin{array}
[c]{c}%
\\
\\
\end{array}
\right.  \\
\\
\\
\\
\\
\\
\\
\end{array}%
\begin{tabular}
[c]{cccccccccccccc}%
${\footnotesize 1}$ & ${\footnotesize 2}$ & ${\footnotesize 3}$ &
${\footnotesize 4}$ & ${\footnotesize 5}$ & ${\footnotesize 6}$ &
${\footnotesize 7}$ & ${\footnotesize 8}$ & ${\footnotesize 9}$ &
${\footnotesize 10}$ & ${\footnotesize 11}$ & ${\footnotesize 12}$ &
${\footnotesize 13}$ & ${\footnotesize 14}$\\\hline
\multicolumn{1}{|c}{$1$} & \multicolumn{1}{|c}{$1$} & \multicolumn{1}{|c}{$1$}
& \multicolumn{1}{|c}{$1$} & \multicolumn{1}{|c}{$1$} &
\multicolumn{1}{|c}{$2$} & \multicolumn{1}{|c}{$2$} & \multicolumn{1}{|c}{$2$}
& \multicolumn{1}{|c}{$2$} & \multicolumn{1}{|c}{$2$} &
\multicolumn{1}{|c}{$3$} & \multicolumn{1}{|c}{$3$} & \multicolumn{1}{|c}{$3$}
& \multicolumn{1}{|c|}{$3$}\\\hline
\multicolumn{1}{|c}{$2$} & \multicolumn{1}{|c}{$3$} & \multicolumn{1}{|c}{$4$}
& \multicolumn{1}{|c}{$5$} & \multicolumn{1}{|c}{$6$} &
\multicolumn{1}{|c}{$4$} & \multicolumn{1}{|c}{$7$} & \multicolumn{1}{|c}{$5$}
& \multicolumn{1}{|c}{$8$} & \multicolumn{1}{|c}{$9$} &
\multicolumn{1}{|c}{$0$} & \multicolumn{1}{|c}{$0$} & \multicolumn{1}{|c}{$0$}
& \multicolumn{1}{|c|}{$9$}\\\hline
${\footnotesize 15}$ & ${\footnotesize 16}$ & ${\footnotesize 17}$ &
${\footnotesize 18}$ & ${\footnotesize 19}$ & ${\footnotesize 20}$ &
${\footnotesize 21}$ & ${\footnotesize 22}$ & ${\footnotesize 23}$ &
${\footnotesize 24}$ & ${\footnotesize 25}$ & ${\footnotesize 26}$ &
${\footnotesize 27}$ & ${\footnotesize 28}$\\\hline
\multicolumn{1}{|c}{$4$} & \multicolumn{1}{|c}{$4$} & \multicolumn{1}{|c}{$4$}
& \multicolumn{1}{|c}{$4$} & \multicolumn{1}{|c}{$5$} &
\multicolumn{1}{|c}{$6$} & \multicolumn{1}{|c}{$6$} & \multicolumn{1}{|c}{$6$}
& \multicolumn{1}{|c}{$6$} & \multicolumn{1}{|c}{$7$} &
\multicolumn{1}{|c}{$7$} & \multicolumn{1}{|c}{$9$} & \multicolumn{1}{|c}{$9$}
& \multicolumn{1}{|c|}{$0$}\\\hline
\multicolumn{1}{|c}{$0$} & \multicolumn{1}{|c}{$6$} & \multicolumn{1}{|c}{$0$}
& \multicolumn{1}{|c}{$0$} & \multicolumn{1}{|c}{$0$} &
\multicolumn{1}{|c}{$8$} & \multicolumn{1}{|c}{$7$} & \multicolumn{1}{|c}{$0$}
& \multicolumn{1}{|c}{$0$} & \multicolumn{1}{|c}{$0$} &
\multicolumn{1}{|c}{$0$} & \multicolumn{1}{|c}{$0$} & \multicolumn{1}{|c}{$0$}
& \multicolumn{1}{|c|}{$0$}\\\hline
${\footnotesize 29}$ & ${\footnotesize \cdots}$ & ${\footnotesize 50}$ &  &  &
&  &  &  &  &  &  &  & \\\cline{1-3}%
\multicolumn{1}{|c}{$0$} & \multicolumn{1}{|c}{$\cdots$} &
\multicolumn{1}{|c}{$0$} & \multicolumn{1}{|c}{} &  &  &  &  &  &  &  &  &  &
\\\cline{1-3}%
\multicolumn{1}{|c}{$0$} & \multicolumn{1}{|c}{$\cdots$} &
\multicolumn{1}{|c}{$0$} & \multicolumn{1}{|c}{} &  &  &  &  &  &  &  &  &  &
\\\cline{1-3}%
\end{tabular}%
\begin{array}
[c]{c}%
\\
\\
\\
\\
\\
\\
\\
\left.
\begin{array}
[c]{c}%
\\
\\
\end{array}
\right]  \\
\end{array}
\label{Eq_construction_example_2}%
\end{equation}%
\begin{figure}
[ptb]
\begin{center}
\includegraphics[
height=3.0562in,
width=2.7899in
]%
{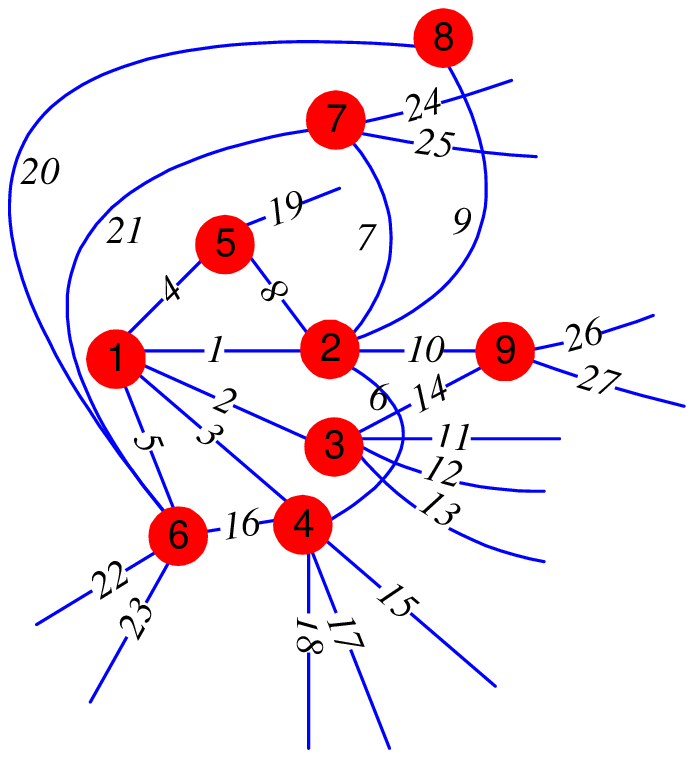}%
\caption{The example of construction. Both or one of the two endnodes of links
$1$-$27$ are determined.}%
\label{2nd_phase}%
\end{center}
\end{figure}
\begin{figure}
[ptb]
\begin{center}
\includegraphics[
height=2.3687in,
width=3.5475in
]%
{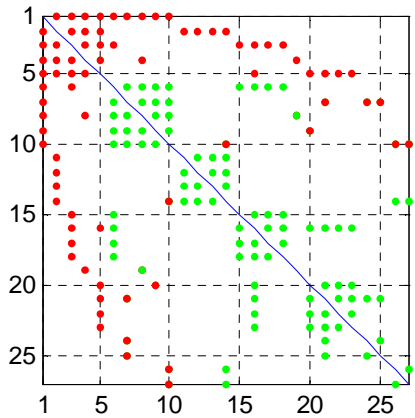}%
\caption{The LAM\ of the constructed part (links $1$-$27)$ of graph are
computed. The green $1$-entries are determined by the red $1$-entries.}%
\label{phase_2nd_matrix}%
\end{center}
\end{figure}

\item Similarly, let node $10$ be the endnode of link $11$, $19$ and $28$%
-$30$, and let node $11$ be the endnode of link $12$, $18$ and $31$-$35$, and
let node $12$ be the endnode of link $13$ and $36$, as shown in Equation
$\left(  \ref{Eq_construction_example_3}\right)  $ and Figure \ref{3rd_phase}
(a).
\begin{equation}%
\begin{array}
[c]{c}%
\\
E=
\end{array}%
\begin{array}
[c]{c}%
\\
\left[
\begin{array}
[c]{c}%
\\
\\
\end{array}
\right.  \\
\\
\\
\\
\\
\\
\\
\end{array}%
\begin{tabular}
[c]{cccccccccccccc}%
${\footnotesize 1}$ & ${\footnotesize 2}$ & ${\footnotesize 3}$ &
${\footnotesize 4}$ & ${\footnotesize 5}$ & ${\footnotesize 6}$ &
${\footnotesize 7}$ & ${\footnotesize 8}$ & ${\footnotesize 9}$ &
${\footnotesize 10}$ & ${\footnotesize 11}$ & ${\footnotesize 12}$ &
${\footnotesize 13}$ & ${\footnotesize 14}$\\\hline
\multicolumn{1}{|c}{$1$} & \multicolumn{1}{|c}{$1$} & \multicolumn{1}{|c}{$1$}
& \multicolumn{1}{|c}{$1$} & \multicolumn{1}{|c}{$1$} &
\multicolumn{1}{|c}{$2$} & \multicolumn{1}{|c}{$2$} & \multicolumn{1}{|c}{$2$}
& \multicolumn{1}{|c}{$2$} & \multicolumn{1}{|c}{$2$} &
\multicolumn{1}{|c}{$3$} & \multicolumn{1}{|c}{$3$} & \multicolumn{1}{|c}{$3$}
& \multicolumn{1}{|c|}{$3$}\\\hline
\multicolumn{1}{|c}{$2$} & \multicolumn{1}{|c}{$3$} & \multicolumn{1}{|c}{$4$}
& \multicolumn{1}{|c}{$5$} & \multicolumn{1}{|c}{$6$} &
\multicolumn{1}{|c}{$4$} & \multicolumn{1}{|c}{$7$} & \multicolumn{1}{|c}{$5$}
& \multicolumn{1}{|c}{$8$} & \multicolumn{1}{|c}{$9$} &
\multicolumn{1}{|c}{$10$} & \multicolumn{1}{|c}{$11$} &
\multicolumn{1}{|c}{$12$} & \multicolumn{1}{|c|}{$9$}\\\hline
${\footnotesize 15}$ & ${\footnotesize 16}$ & ${\footnotesize 17}$ &
${\footnotesize 18}$ & ${\footnotesize 19}$ & ${\footnotesize 20}$ &
${\footnotesize 21}$ & ${\footnotesize 22}$ & ${\footnotesize 23}$ &
${\footnotesize 24}$ & ${\footnotesize 25}$ & ${\footnotesize 26}$ &
${\footnotesize 27}$ & ${\footnotesize 28}$\\\hline
\multicolumn{1}{|c}{$4$} & \multicolumn{1}{|c}{$4$} & \multicolumn{1}{|c}{$4$}
& \multicolumn{1}{|c}{$4$} & \multicolumn{1}{|c}{$5$} &
\multicolumn{1}{|c}{$6$} & \multicolumn{1}{|c}{$6$} & \multicolumn{1}{|c}{$6$}
& \multicolumn{1}{|c}{$6$} & \multicolumn{1}{|c}{$7$} &
\multicolumn{1}{|c}{$7$} & \multicolumn{1}{|c}{$9$} & \multicolumn{1}{|c}{$9$}
& \multicolumn{1}{|c|}{$10$}\\\hline
\multicolumn{1}{|c}{$0$} & \multicolumn{1}{|c}{$6$} & \multicolumn{1}{|c}{$0$}
& \multicolumn{1}{|c}{$11$} & \multicolumn{1}{|c}{$10$} &
\multicolumn{1}{|c}{$8$} & \multicolumn{1}{|c}{$7$} & \multicolumn{1}{|c}{$0$}
& \multicolumn{1}{|c}{$0$} & \multicolumn{1}{|c}{$0$} &
\multicolumn{1}{|c}{$0$} & \multicolumn{1}{|c}{$0$} & \multicolumn{1}{|c}{$0$}
& \multicolumn{1}{|c|}{$0$}\\\hline
${\footnotesize 29}$ & ${\footnotesize 30}$ & ${\footnotesize 31}$ &
${\footnotesize 32}$ & ${\footnotesize 33}$ & ${\footnotesize 34}$ &
${\footnotesize 35}$ & ${\footnotesize 36}$ & ${\footnotesize 37}$ &
${\footnotesize \cdots}$ & ${\footnotesize 50}$ &  &  & \\\cline{1-11}%
\multicolumn{1}{|c}{$10$} & \multicolumn{1}{|c}{$10$} &
\multicolumn{1}{|c}{$11$} & \multicolumn{1}{|c}{$11$} &
\multicolumn{1}{|c}{$11$} & \multicolumn{1}{|c}{$11$} &
\multicolumn{1}{|c}{$11$} & \multicolumn{1}{|c}{$12$} &
\multicolumn{1}{|c}{$0$} & \multicolumn{1}{|c}{$\cdots$} &
\multicolumn{1}{|c}{$0$} & \multicolumn{1}{|c}{} &  & \\\cline{1-11}%
\multicolumn{1}{|c}{$0$} & \multicolumn{1}{|c}{$0$} & \multicolumn{1}{|c}{$0$}
& \multicolumn{1}{|c}{$0$} & \multicolumn{1}{|c}{$0$} &
\multicolumn{1}{|c}{$0$} & \multicolumn{1}{|c}{$0$} & \multicolumn{1}{|c}{$0$}
& \multicolumn{1}{|c}{$0$} & \multicolumn{1}{|c}{$\cdots$} &
\multicolumn{1}{|c}{$0$} & \multicolumn{1}{|c}{} &  & \\\cline{1-11}%
\end{tabular}%
\begin{array}
[c]{c}%
\\
\\
\\
\\
\\
\\
\\
\left.
\begin{array}
[c]{c}%
\\
\\
\end{array}
\right]  \\
\end{array}
\label{Eq_construction_example_3}%
\end{equation}%
\begin{figure}
[ptb]
\begin{center}
\includegraphics[
height=3.1981in,
width=3.026in
]%
{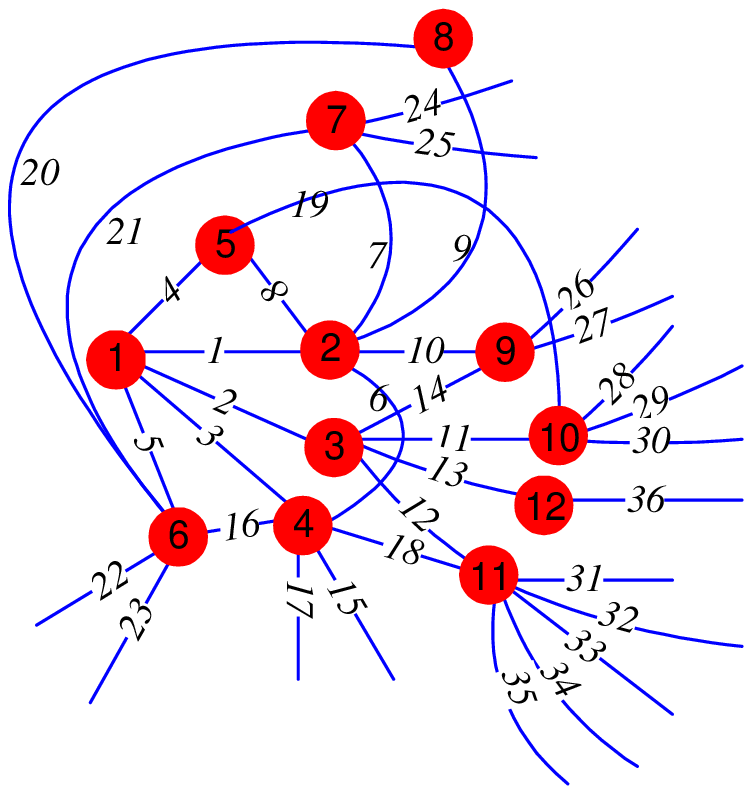}%
\caption{The example of construction. Both or one of the two endnodes of links
$1$-$36$ are determined.}%
\label{3rd_phase}%
\end{center}
\end{figure}
\begin{figure}
[ptb]
\begin{center}
\includegraphics[
height=3.0813in,
width=3.3797in
]%
{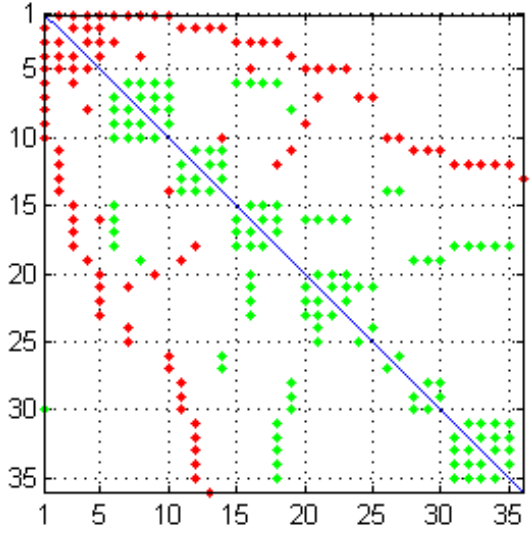}%
\caption{The LAM\ of the constructed part (links $1$-$36)$ of graph are
computed. The green $1$-entries are determined by the red $1$-entries.}%
\label{phase_3rd_matrix}%
\end{center}
\end{figure}

\item Constructing in this way, the two endnodes of all the links are
eventually determined, as shown in Equation $\left(
\ref{Eq_construction_example_done}\right)  $ and Figure \ref{4th_phase} (a).
The final structure of the matrix $E$ exhibits the link list of the original
graph $G$ which consists of $30$ nodes and $50$ links. For example, link $36$
connects node $12$ and node $15$ in $G$. The matrix $E$ is readily transformed
into the adjacency matrix of $G$.%
\begin{equation}%
\begin{array}
[c]{c}%
\\
E=
\end{array}%
\begin{array}
[c]{c}%
\\
\left[
\begin{array}
[c]{c}%
\\
\\
\end{array}
\right.  \\
\\
\\
\\
\\
\\
\\
\\
\\
\\
\end{array}%
\begin{tabular}
[c]{cccccccccccccc}%
${\footnotesize 1}$ & ${\footnotesize 2}$ & ${\footnotesize 3}$ &
${\footnotesize 4}$ & ${\footnotesize 5}$ & ${\footnotesize 6}$ &
${\footnotesize 7}$ & ${\footnotesize 8}$ & ${\footnotesize 9}$ &
${\footnotesize 10}$ & ${\footnotesize 11}$ & ${\footnotesize 12}$ &
${\footnotesize 13}$ & ${\footnotesize 14}$\\\hline
\multicolumn{1}{|c}{$1$} & \multicolumn{1}{|c}{$1$} & \multicolumn{1}{|c}{$1$}
& \multicolumn{1}{|c}{$1$} & \multicolumn{1}{|c}{$1$} &
\multicolumn{1}{|c}{$2$} & \multicolumn{1}{|c}{$2$} & \multicolumn{1}{|c}{$2$}
& \multicolumn{1}{|c}{$2$} & \multicolumn{1}{|c}{$2$} &
\multicolumn{1}{|c}{$3$} & \multicolumn{1}{|c}{$3$} & \multicolumn{1}{|c}{$3$}
& \multicolumn{1}{|c|}{$3$}\\\hline
\multicolumn{1}{|c}{$2$} & \multicolumn{1}{|c}{$3$} & \multicolumn{1}{|c}{$4$}
& \multicolumn{1}{|c}{$5$} & \multicolumn{1}{|c}{$6$} &
\multicolumn{1}{|c}{$4$} & \multicolumn{1}{|c}{$7$} & \multicolumn{1}{|c}{$5$}
& \multicolumn{1}{|c}{$8$} & \multicolumn{1}{|c}{$9$} &
\multicolumn{1}{|c}{$10$} & \multicolumn{1}{|c}{$11$} &
\multicolumn{1}{|c}{$12$} & \multicolumn{1}{|c|}{$9$}\\\hline
${\footnotesize 15}$ & ${\footnotesize 16}$ & ${\footnotesize 17}$ &
${\footnotesize 18}$ & ${\footnotesize 19}$ & ${\footnotesize 20}$ &
${\footnotesize 21}$ & ${\footnotesize 22}$ & ${\footnotesize 23}$ &
${\footnotesize 24}$ & ${\footnotesize 25}$ & ${\footnotesize 26}$ &
${\footnotesize 27}$ & ${\footnotesize 28}$\\\hline
\multicolumn{1}{|c}{$4$} & \multicolumn{1}{|c}{$4$} & \multicolumn{1}{|c}{$4$}
& \multicolumn{1}{|c}{$4$} & \multicolumn{1}{|c}{$5$} &
\multicolumn{1}{|c}{$6$} & \multicolumn{1}{|c}{$6$} & \multicolumn{1}{|c}{$6$}
& \multicolumn{1}{|c}{$6$} & \multicolumn{1}{|c}{$7$} &
\multicolumn{1}{|c}{$7$} & \multicolumn{1}{|c}{$9$} & \multicolumn{1}{|c}{$9$}
& \multicolumn{1}{|c|}{$10$}\\\hline
\multicolumn{1}{|c}{$13$} & \multicolumn{1}{|c}{$6$} &
\multicolumn{1}{|c}{$14$} & \multicolumn{1}{|c}{$11$} &
\multicolumn{1}{|c}{$10$} & \multicolumn{1}{|c}{$8$} & \multicolumn{1}{|c}{$7$%
} & \multicolumn{1}{|c}{$15$} & \multicolumn{1}{|c}{$16$} &
\multicolumn{1}{|c}{$13$} & \multicolumn{1}{|c}{$17$} &
\multicolumn{1}{|c}{$18$} & \multicolumn{1}{|c}{$19$} &
\multicolumn{1}{|c|}{$19$}\\\hline
${\footnotesize 29}$ & ${\footnotesize 30}$ & ${\footnotesize 31}$ &
${\footnotesize 32}$ & ${\footnotesize 33}$ & ${\footnotesize 34}$ &
${\footnotesize 35}$ & ${\footnotesize 36}$ & ${\footnotesize 37}$ &
${\footnotesize 38}$ & ${\footnotesize 39}$ & ${\footnotesize 40}$ &
${\footnotesize 41}$ & ${\footnotesize 42}$\\\hline
\multicolumn{1}{|c}{$10$} & \multicolumn{1}{|c}{$10$} &
\multicolumn{1}{|c}{$11$} & \multicolumn{1}{|c}{$11$} &
\multicolumn{1}{|c}{$11$} & \multicolumn{1}{|c}{$11$} &
\multicolumn{1}{|c}{$11$} & \multicolumn{1}{|c}{$12$} &
\multicolumn{1}{|c}{$13$} & \multicolumn{1}{|c}{$14$} &
\multicolumn{1}{|c}{$15$} & \multicolumn{1}{|c}{$15$} &
\multicolumn{1}{|c}{$16$} & \multicolumn{1}{|c|}{$17$}\\\hline
\multicolumn{1}{|c}{$20$} & \multicolumn{1}{|c}{$21$} &
\multicolumn{1}{|c}{$22$} & \multicolumn{1}{|c}{$23$} &
\multicolumn{1}{|c}{$24$} & \multicolumn{1}{|c}{$21$} &
\multicolumn{1}{|c}{$16$} & \multicolumn{1}{|c}{$25$} &
\multicolumn{1}{|c}{$23$} & \multicolumn{1}{|c}{$26$} &
\multicolumn{1}{|c}{$17$} & \multicolumn{1}{|c}{$26$} &
\multicolumn{1}{|c}{$24$} & \multicolumn{1}{|c|}{$24$}\\\hline
${\footnotesize 43}$ & ${\footnotesize 44}$ & ${\footnotesize 45}$ &
${\footnotesize 46}$ & ${\footnotesize 47}$ & ${\footnotesize 48}$ &
${\footnotesize 49}$ & ${\footnotesize 50}$ &  &  &  &  &  & \\\cline{1-8}%
\multicolumn{1}{|c}{$17$} & \multicolumn{1}{|c}{$17$} &
\multicolumn{1}{|c}{$18$} & \multicolumn{1}{|c}{$23$} &
\multicolumn{1}{|c}{$23$} & \multicolumn{1}{|c}{$23$} &
\multicolumn{1}{|c}{$27$} & \multicolumn{1}{|c}{$27$} & \multicolumn{1}{|c}{}
&  &  &  &  & \\\cline{1-8}%
\multicolumn{1}{|c}{$27$} & \multicolumn{1}{|c}{$28$} &
\multicolumn{1}{|c}{$21$} & \multicolumn{1}{|c}{$26$} &
\multicolumn{1}{|c}{$25$} & \multicolumn{1}{|c}{$28$} &
\multicolumn{1}{|c}{$29$} & \multicolumn{1}{|c}{$30$} & \multicolumn{1}{|c}{}
&  &  &  &  & \\\cline{1-8}%
\end{tabular}%
\begin{array}
[c]{c}%
\\
\\
\\
\\
\\
\\
\\
\\
\\
\\
\left.
\begin{array}
[c]{c}%
\\
\\
\end{array}
\right]  \\
\end{array}
\label{Eq_construction_example_done}%
\end{equation}%
\begin{figure}
[ptb]
\begin{center}
\includegraphics[
height=3.4601in,
width=3.6037in
]%
{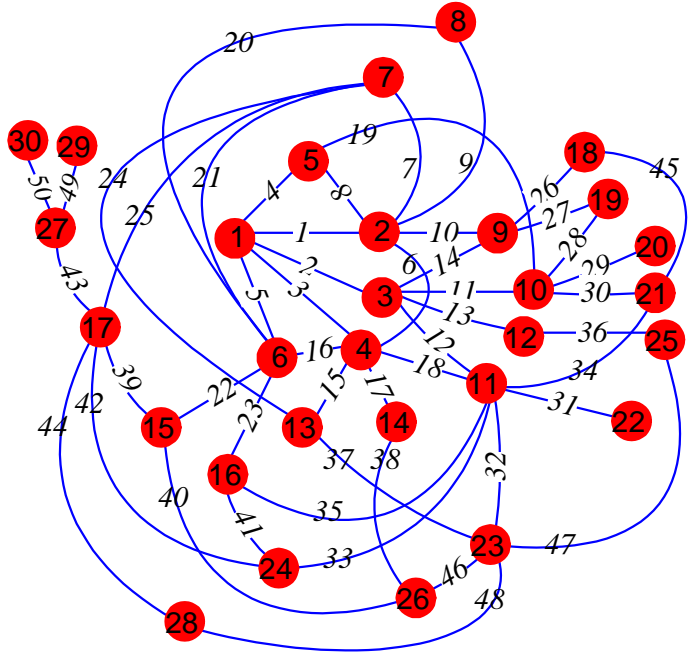}%
\caption{The example of construction. The two endnodes of all links are
determined.}%
\label{4th_phase}%
\end{center}
\end{figure}
\begin{figure}
[ptb]
\begin{center}
\includegraphics[
height=3.6443in,
width=3.6867in
]%
{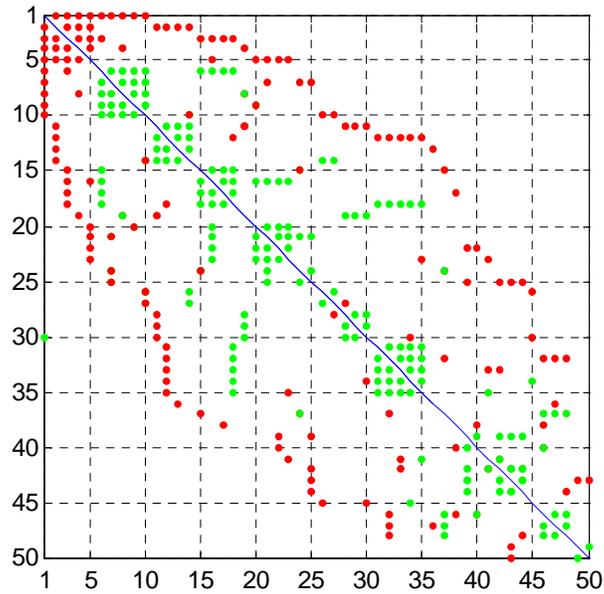}%
\caption{The LAM\ of the constructed graph is computed. The green $1$-entries
are determined by the red $1$-entries.}%
\label{phase_last_matrix}%
\end{center}
\end{figure}

\end{enumerate}

\subsubsection{Initialization (The recognition of the endnodes of the first
link and its neighboring links)}

\label{Section_construction_algo_intial}%

\begin{algorithm}%

\caption{$E_{2 \times L_{G}} \Leftarrow Initialization(C,s_{1},s_{2},s_{3})$}\label{algo_initialization}%

\begin{algorithmic}[1]
\IF{$s_{3} \geq 1$}
\STATE $E \Leftarrow \mathcal{E}$
\ELSIF{$s_{1}=1$}
\STATE $E \Leftarrow \mathcal{E}_{1}$
\ELSIF{$s_{1}=2$}
\STATE $E_{2 \times L_{G}} \Leftarrow Initialization2(C,s_{2},s_{3})$
\ELSIF{$s_{1}=3$}
\STATE $E_{2 \times L_{G}} \Leftarrow Initialization3(C,s_{2},s_{3})$
\ELSIF{$s_{1} \geq 4$}
\STATE $E_{2 \times L_{G}} \Leftarrow Initialization4(C,s_{1},s_{2},s_{3})$
\ENDIF
\end{algorithmic}%

\end{algorithm}%

When $s_{3}\geq1$, Theorem \ref{theorem_relabeling_three_lines} implies that
$s_{2}\geq2$, $s_{1}\geq3$ and links $2,3,\cdots,s_{3}+3$ are incident to the
left endnode of link $1$ and links $s_{3}+4,\cdots,s_{1}+1$ are incident to
the right endnode of link $1$. Therefore, line 1-2 of Algorithm
\ref{algo_initialization} initialize $E$ by $\mathcal{E}$. The numbers above
the matrix $\mathcal{E}$ in $\left(  \ref{Eq_E5}\right)  $ are the column
numbers, which indicate the link numbers, and $\mathcal{E}$\ has the following
structure,
\begin{equation}%
\begin{array}
[c]{c}%
\\
\mathcal{E}=
\end{array}%
\begin{array}
[c]{c}%
\\
\left[
\begin{array}
[c]{c}%
\\
\\
\end{array}
\right.
\end{array}%
\begin{tabular}
[c]{cccccccccc}%
${\footnotesize 1}$ & ${\footnotesize 2}$ & ${\footnotesize \cdots}$ &
${\footnotesize s}_{3}{\footnotesize +3}$ & ${\footnotesize s}_{3}%
{\footnotesize +4}$ & ${\footnotesize \cdots}$ & ${\footnotesize s}%
_{1}{\footnotesize +1}$ & ${\footnotesize s}_{1}{\footnotesize +2}$ &
${\footnotesize \cdots}$ & {\footnotesize $L_{G}$}\\\hline
\multicolumn{1}{|c}{$1$} & \multicolumn{1}{|c}{$1$} &
\multicolumn{1}{|c}{$\cdots$} & \multicolumn{1}{|c}{$1$} &
\multicolumn{1}{|c}{$2$} & \multicolumn{1}{|c}{$\cdots$} &
\multicolumn{1}{|c}{$2$} & \multicolumn{1}{|c}{$0$} &
\multicolumn{1}{|c}{$\cdots$} & \multicolumn{1}{|c|}{$0$}\\\hline
\multicolumn{1}{|c}{$2$} & \multicolumn{1}{|c}{$0$} &
\multicolumn{1}{|c}{$\cdots$} & \multicolumn{1}{|c}{$0$} &
\multicolumn{1}{|c}{$0$} & \multicolumn{1}{|c}{$\cdots$} &
\multicolumn{1}{|c}{$0$} & \multicolumn{1}{|c}{$0$} &
\multicolumn{1}{|c}{$\cdots$} & \multicolumn{1}{|c|}{$0$}\\\hline
\end{tabular}%
\begin{array}
[c]{c}%
\\
\left.
\begin{array}
[c]{c}%
\\
\\
\end{array}
\right]
\end{array}
\label{Eq_E5}%
\end{equation}

When $s_{3}=0$, Theorem \ref{theorem_relabeling_three_lines} cannot be used.
However, the limited number of cases of $s_{3}=0$ enables us to accomplish the
initialization with the detailed analysis in the Appendix
\ref{Section_appendix_A}.

\subsubsection{The recognition of the endnodes of the whole graph}

\label{Section_construction_algo_}Lines 1-2 of Algorithm \ref{algo_Marinlinga}
relabel the given LAM $C$ and determine the initial state. In the initial
state, link $1$ is always incident to node $1$ and $2$. Some of the
neighboring links of link $1$ are incident to node $1$, and the other
neighboring links are incident to node $2$. The second endnodes of the
neighboring links of link $1$ have not decided yet in the initial state.

Line 3 initiates the number of nodes $N_{G}$ to $2$. The two endnodes of link
$1$ are already determined. Starting with link $2$ until link $L_{G}$ (line
4), the number of nodes $N_{G}$ increases by $1$ (line 6) if the second
endnode of link $i$ is not determined (line 5). Let the second endnode of link
$i$ be $N_{G}$ (line 7). When link $i$ is adjacent to link $j$, $j=i+1,\cdots
,L_{G}$ (lines 8-9), let the first endnode of link $j$ be $N_{G}$ (line 11) if
the first endnode of link $j$ is not determined (line 10). If the first
endnode of link $j$ is determined but the second endnode is not determined and
links $i$ and $j$ do not share the first endnode (line 12), let the second
endnode of link $j$ be $N_{G}$ (line 13).%

\begin{algorithm}%

\caption{$E_{2 \times L_{G}} \Leftarrow MARINLINGA(C)$}\label{algo_Marinlinga}

\begin{algorithmic}[1]
\STATE $(C,s_{1},s_{2},s_{3}) \Leftarrow MatrixRelabeling(C)$
\STATE $E_{2 \times L_{G}} \Leftarrow Initialization(C,s_{1},s_{2},s_{3})$
\STATE $N \Leftarrow 2$
\FOR{$i=2$ to $L_{G}$}
\IF{$e_{2i} = 0$}
\STATE $N \Leftarrow N+1$
\STATE $e_{2i} \Leftarrow N$
\FOR{$j=i+1$ to $L_{G}$}
\IF{$c_{ij} = 1$}
\IF{$e_{1j} = 0$}
\STATE $e_{1j} \Leftarrow N$
\ELSIF{$e_{2j} = 0$ and $e_{1i} \neq e_{1j}$}
\STATE $e_{2j} \Leftarrow N$
\ENDIF
\ENDIF
\ENDFOR
\ENDIF
\ENDFOR
\end{algorithmic}%

\end{algorithm}%

\subsection{Worst case complexity of MARINLINGA}

Algorithm \ref{algo_label_swapping} has a complexity of $O\left(
L_{G}\right)  $. The complexity of Algorithm \ref{algo_link_relabeling} can be
computed as follows. Line 1 has a complexity of $O\left(  L_{G}\right)  $. In
the worst case, the function of line 2, Algorithm
\ref{algo_group_label_swapping} has a complexity of $O\left(  L_{G}%
^{2}\right)  $, if $m$ in line 15 of Algorithm \ref{algo_group_label_swapping}
is proportional to $L_{G}$. The worst case complexity of lines 3-6 is also
$O\left(  L_{G}^{2}\right)  $. Hence, lines 1-6 have a complexity of $O\left(
L_{G}^{2}\right)  $. Neglect $O\left(  1\right)  $ operations of lines 7-8.
The times that lines 9-14 are executed is stored in $k$. If $k$ is
proportional to $L_{G}$, $m$ in line 15 of Algorithm
\ref{algo_group_label_swapping} must be bounded by a constant, then the
complexity of line 11 is $O\left(  L_{G}\right)  $. If $k$ is bounded, the
complexity of line 11 will be $O\left(  L_{G}^{2}\right)  $. Therefore, lines
9-14 have a worst case complexity of $O\left(  L_{G}^{2}\right)  $. Hence, the
complexity of Algorithm \ref{algo_link_relabeling} is $O\left(  L_{G}%
^{2}\right)  $.

Algorithm \ref{algo_initialization_s1_2}, \ref{algo_initialization_s1_3} and
\ref{algo_initialization_s1_4} have a worst case complexity of $O\left(
1\right)  $, hence the complexity of Algorithm \ref{algo_initialization} is
also $O\left(  1\right)  $. Lines 4-18 of the main Algorithm
\ref{algo_Marinlinga} have a worst case complexity of $O\left(  L_{G}%
^{2}\right)  $. In summary, the worst case complexity of the MARINLINGA is
$O\left(  L_{G}^{2}\right)  $. Since the number of links of \emph{the original
graph} $G$\ and the number of nodes of \emph{the line graph} $l\left(
G\right)  $\ are equal, $L_{G}=N_{l\left(  G\right)  }$, the worst case
complexity of the MARINLINGA is written as $O\left(  N_{l\left(  G\right)
}^{2}\right)  $.

\section{Comparison with Roussopoulos's algorithm}

\begin{figure}
[ptb]
\begin{center}
\includegraphics[
height=5.3039in,
width=5.5832in
]%
{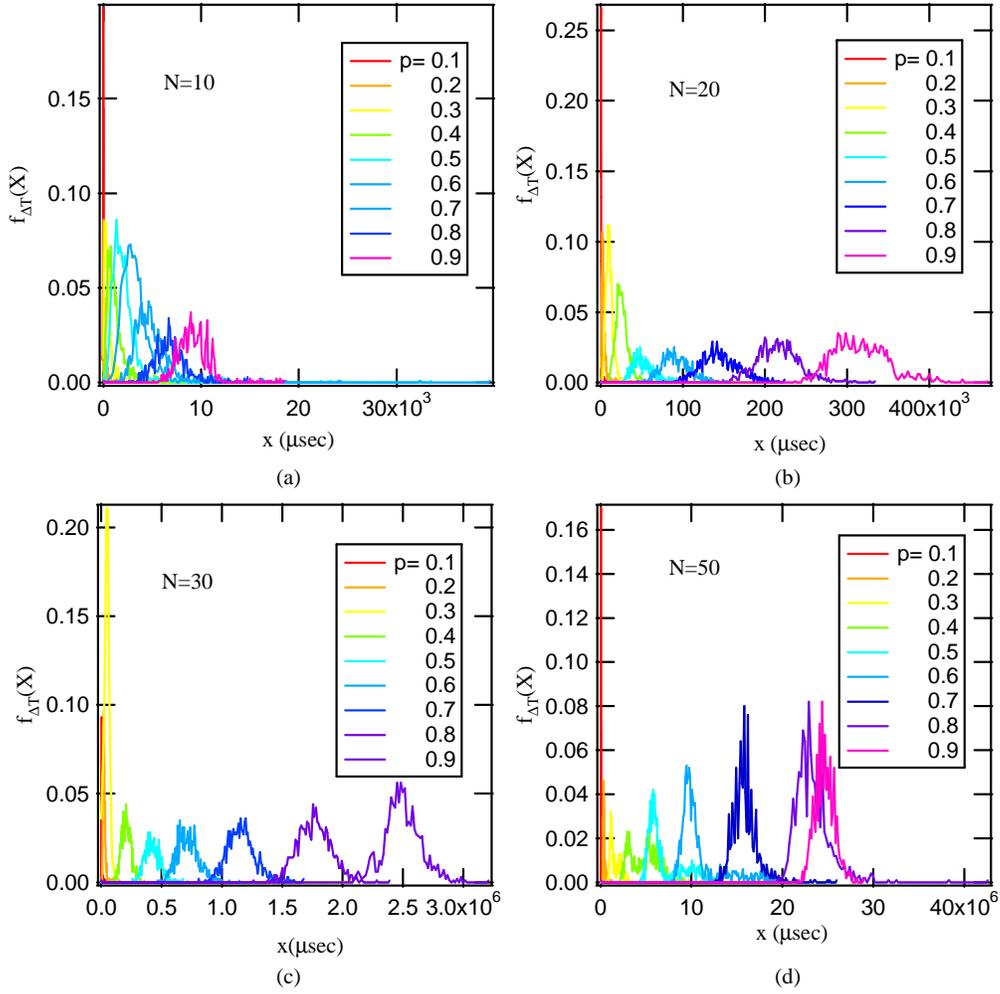}%
\caption{PDFs of the $\Delta T$ for $N=10,20,30$ and $50$ and
$p=\{0.1,0.2,\ldots,0.9\}$}%
\label{n10n20n30n50}%
\end{center}
\end{figure}

We use the same input line graphs for both MARINLINGA and Roussopoulos's
algorithm. We start with line graphs constructed from Erd\H{o}s-R\'{e}nyi
random graphs $G_{p}\left(  N_{G}\right)  $ \cite{ER}. We calculate the
probability density function of the difference $\Delta T$ between the running
time of Roussopoulos's algorithm ($T_{\text{Roussopoulos}}$) and MARINLINGA
($T_{\text{MARINLINGA}}$)%

\begin{equation}
\Delta T=T_{\text{Roussopoulos}}-T_{\text{MARINLINGA}} \label{time_difference}%
\end{equation}
We randomly create $1000$ different line graphs based on the class of
Erd\"{o}s-R\'{e}nyi random graphs $G_{p}\left(  N_{G}\right)  $ for each
number of nodes $N_{G}=\{10,20,30,50\}$ and link density
$p=\{0.1,0.2,...,0.9\}$. The probability density functions of the time
difference $f_{\Delta T}\left(  x\right)  $ for each class of line graphs of
$G_{p}\left(  N_{G}\right)  $ are shown in Figure \ref{n10n20n30n50}.%
The values of the probability density function are nearly always positive
which means in practice that MARINLINGA needs less time for the execution than
Roussopoulos's algorithm.

We calculate the expectation according to \cite{PVM_PerformanceAnalysisCUP}
and the experimental results for all of the mentioned cases $\Delta T_{i}$ for
$i=1,2,\ldots,1000$.
\begin{equation}
E[\Delta T]=\sum k\text{ Pr}[\Delta T=k]=\frac{\sum_{i=1}^{1000}\Delta T_{i}%
}{1000} \label{expectation_time_difference}%
\end{equation}
The results in milliseconds are given in the Table
\ref{table_time_diff_expectation}.%

\begin{table}[tbp] \centering
\begin{tabular}
[c]{||lc|c|c|c|c||}\hline\hline
& $N_{G}$ & $10$ & $20$ & $30$ & $50$\\
$\ \ p$ &  &  &  &  & \\\hline\hline
\multicolumn{2}{||c}{$.1$} & \multicolumn{1}{|r|}{$0.0462$} &
\multicolumn{1}{|r|}{$0.7839$} & \multicolumn{1}{|r|}{$3.6660$} &
\multicolumn{1}{|r||}{$42.8153$}\\\hline
\multicolumn{2}{||c}{$.2$} & \multicolumn{1}{|r|}{$0.1682$} &
\multicolumn{1}{|r|}{$4.3012$} & \multicolumn{1}{|r|}{$26.1086$} &
\multicolumn{1}{|r||}{$385.3208$}\\\hline
\multicolumn{2}{||c}{$.3$} & \multicolumn{1}{|r|}{$0.6281$} &
\multicolumn{1}{|r|}{$11.9492$} & \multicolumn{1}{|r|}{$49.826$} &
\multicolumn{1}{|r||}{$1582.5827$}\\\hline
\multicolumn{2}{||c}{$.4$} & \multicolumn{1}{|r|}{$1.3921$} &
\multicolumn{1}{|r|}{$26.4580$} & \multicolumn{1}{|r|}{$209.3056$} &
\multicolumn{1}{|r||}{$4802.8110$}\\\hline
\multicolumn{2}{||c}{$.5$} & \multicolumn{1}{|r|}{$2.4269$} &
\multicolumn{1}{|r|}{$51.8393$} & \multicolumn{1}{|r|}{$422.8641$} &
\multicolumn{1}{|r||}{$7060.6021$}\\\hline
\multicolumn{2}{||c}{$.6$} & \multicolumn{1}{|r|}{$3.5808$} &
\multicolumn{1}{|r|}{$91.8258$} & \multicolumn{1}{|r|}{$719.6185$} &
\multicolumn{1}{|r||}{$10978.1532$}\\\hline
\multicolumn{2}{||c}{$.7$} & \multicolumn{1}{|r|}{$4.9181$} &
\multicolumn{1}{|r|}{$145.5864$} & \multicolumn{1}{|r|}{$1165.1978$} &
\multicolumn{1}{|r||}{$15942.7395$}\\\hline
\multicolumn{2}{||c}{$.8$} & \multicolumn{1}{|r|}{$6.8952$} &
\multicolumn{1}{|r|}{$217.2041$} & \multicolumn{1}{|r|}{$1768.9901$} &
\multicolumn{1}{|r||}{$23393.5148$}\\\hline
\multicolumn{2}{||c}{$.9$} & \multicolumn{1}{|r|}{$9.2986$} &
\multicolumn{1}{|r|}{$317.4966$} & \multicolumn{1}{|r|}{$2501.4799$} &
\multicolumn{1}{|r||}{$24723.4739$}\\\hline\hline
\end{tabular}
\caption{Expectations of the time difference ($\mu$sec)}\label{table_time_diff_expectation}%
\end{table}%

Additionally, we calculate the probability that MARINLINGA is slower than
Roussopoulos's algorithm: Pr$[\Delta T<0]$ for each $N_{G}$ and $p$. The
simulation shows that Pr$[\Delta T<0]>0$ only for $N_{G}=10$ and $p\leq0.2$,
in which the graphs are mostly disconnected. When the graph is disconnected,
MARINLINGA needs extra time to partition the graphs into connected components
(see footnote 3). For $N_{G}=10$ and $p=0.1$, Pr$[\Delta T<0]=0.33$ and for
$N_{G}=10$ and $p=0.2$, Pr$[\Delta T<0]=0.01$. For all the other cases%

\begin{equation}
\text{Pr}[\Delta T<0]<0.001
\end{equation}
which means that MARINLINGA is generally more efficient than Roussopoulos's
algorithm. The algorithm to find the maximal connected common subgraphs in
graphs is frequently used in the Roussopoulos's algorithm. This algorithm
requires a high running time, because the problem of finding the maximal
connected common subgraphs is $NP$-complete \cite{VISMARA}. The dependence on
this $NP$-complete algorithm is most significant weakness of Roussopoulos's algorithm.

\section{Conclusion}

We have presented a new algorithm MARINLINGA for reverse line graph
construction. By introducing the concept of LAM, we transformed the problem of
reverse line graph construction into the problem of constructing a graph from
the LAM. MARINLINGA consists of two sub-algorithms: the matrix relabeling
algorithm and the construction algorithm. The matrix relabeling algorithm
preprocesses the LAM into the special order by which we can determine the
neighboring links of the first link and the endnodes of the first link
incident to the neighboring links. The construction algorithm makes the first
two nodes be the endnodes of the first link by default, and thereafter,
determines the endnodes of the remaining links. MARINLINGA has a worst case
complexity of $O(N_{l\left(  G\right)  }^{2})$, where $N_{l\left(  G\right)
}$ denotes the number of nodes of the line graph. We have demonstrated that
MARINLINGA is more time-efficient compared to Roussopoulos's algorithm for
connected line graphs. The complexity of Roussopoulos's algorithm mentioned in
\cite{Roussopoulos} is $O(N_{l\left(  G\right)  }+L_{l\left(  G\right)  })$,
where $N_{l\left(  G\right)  }$ and $L_{l\left(  G\right)  }$ are number of
nodes and links of the line graph. Since $L_{l\left(  G\right)  }=O\left(
N_{l\left(  G\right)  }^{2}\right)  $ in worst case, the complexity of
Roussopoulos's algorithm is also $O(N_{l\left(  G\right)  }^{2})$ in worst
case. However, this analysis neglects the computational time of a
sub-algorithm that determines the maximal connected common subgraph in each
iteration. Finding a maximally connected common subgraph is an $NP$-complete
problem \cite{VISMARA}, implying that Roussopoulos's algorithm is, in fact,
not polynomial in worst case.

\bibliographystyle{plain}
\bibliography{Dajie_linegraph_recons}

\appendix{}

\section{The initialization of the construction algorithm when $s_{3}=0$}

\label{Section_appendix_A}Theorem \ref{theorem_relabeling_three_lines} cannot
be used when $s_{3}=0$. Since there exists limited number of cases of
$s_{3}=0$, we can still accomplish the initialization.

\subsection{When $s_{1}=1$}

Link $1$ has only one right neighboring link: link $2$. Link $1$ does not have
left neighboring links. The initial state of $E$ is $\mathcal{E}_{1}$. Lines
3-4 of Algorithm \ref{algo_initialization} initialize $E$ by $\mathcal{E}_{1}%
$.%
\begin{equation}
\mathcal{E}_{1}=\left[
\begin{tabular}
[c]{ccccc}%
$1$ & $2$ & $0$ & $\cdots$ & $0$\\
$2$ & $0$ & $0$ & $\cdots$ & $0$%
\end{tabular}
\ \ \ \ \ \right]
\end{equation}

\subsection{When $s_{1}=2$}

\begin{figure}
[ptb]
\begin{center}
\includegraphics[
height=2.6091in,
width=2.7674in
]%
{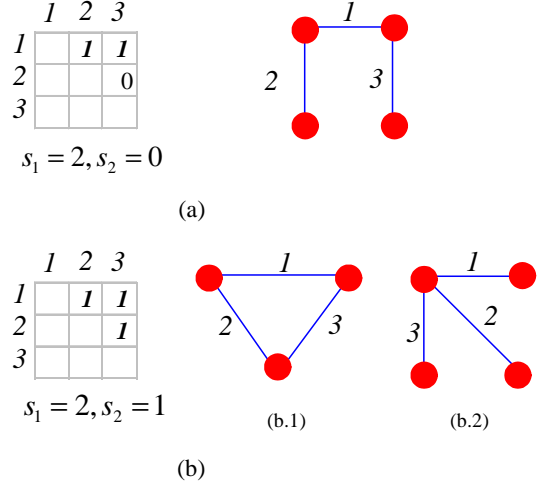}%
\caption{The adjacency patterns of link $1$ and its neighboring links when
$s_{1}=2$. The graphs on the right are the possible configurations
correspondingly.}%
\label{s1smaller4_2}%
\end{center}
\end{figure}

There are different adjacency patterns. The submatrix of $C$ in Figure
\ref{s1smaller4_2} (a) implies that, links $2$ and $3$ are adjacent to link
$1$, and link $2$ is not not adjacent to link $3$. Links $2$ and $3$ must be
incident to two different endnodes of link $1$. The pattern in Figure
\ref{s1smaller4_2} (b) has two possible configurations $K_{3}$ and $K_{1,3}$.
If $s_{1}=2$ and $s_{2}=0$, the initial state is $\mathcal{E}_{2,a}$, as shown
in lines 1-2 of Algorithm \ref{algo_initialization_s1_2}. When $s_{1}=2$ and
$s_{2}=1$, because the graph is connected, either $c_{2,4}=c_{3,4}=1$ or
$c_{2,4}=0,c_{3,4}=1$ or $c_{2,4}=1,c_{3,4}=0$. If $c_{2,4}=c_{3,4}=1$, the
initial state is $\mathcal{E}_{2,b.1}$, which is $K_{3}$, otherwise the
initial state is $\mathcal{E}_{2,b.2}$, which is $K_{1,3}$, as shown in lines
8-12 of Algorithm \ref{algo_initialization_s1_2}.%
\begin{equation}
\mathcal{E}_{2,a}=\left[
\begin{tabular}
[c]{cccccc}%
$1$ & $1$ & $2$ & $0$ & $\cdots$ & $0$\\
$2$ & $0$ & $0$ & $0$ & $\cdots$ & $0$%
\end{tabular}
\ \ \ \ \ \right]
\end{equation}%
\begin{equation}
\mathcal{E}_{2,b.1}=\left[
\begin{tabular}
[c]{cccccc}%
$1$ & $1$ & $2$ & $0$ & $\cdots$ & $0$\\
$2$ & $0$ & $0$ & $0$ & $\cdots$ & $0$%
\end{tabular}
\ \ \ \ \ \right]  ,\mathcal{E}_{2,b.2}=\left[
\begin{tabular}
[c]{cccccc}%
$1$ & $1$ & $1$ & $0$ & $\cdots$ & $0$\\
$2$ & $0$ & $0$ & $0$ & $\cdots$ & $0$%
\end{tabular}
\ \ \ \ \ \right]
\end{equation}
%

\begin{algorithm}%

\caption{$E_{2 \times L_{G}} \Leftarrow Initialization2(C,s_{2})$}\label{algo_initialization_s1_2}%

\begin{algorithmic}[1]
\IF{$s_{2}=0$}
\STATE $E \Leftarrow \mathcal{E}_{2,a}$
\ELSE
\IF{$c_{2,4}=1$ and $c_{3,4}=1$}
\STATE $E \Leftarrow \mathcal{E}_{2,b.2}$
\ELSE
\STATE $E \Leftarrow \mathcal{E}_{2,b.1}$
\ENDIF
\ENDIF
\end{algorithmic}%

\end{algorithm}%

\subsection{When $s_{1}=3$}

\begin{figure}
[ptb]
\begin{center}
\includegraphics[
height=3.32in,
width=5.412in
]%
{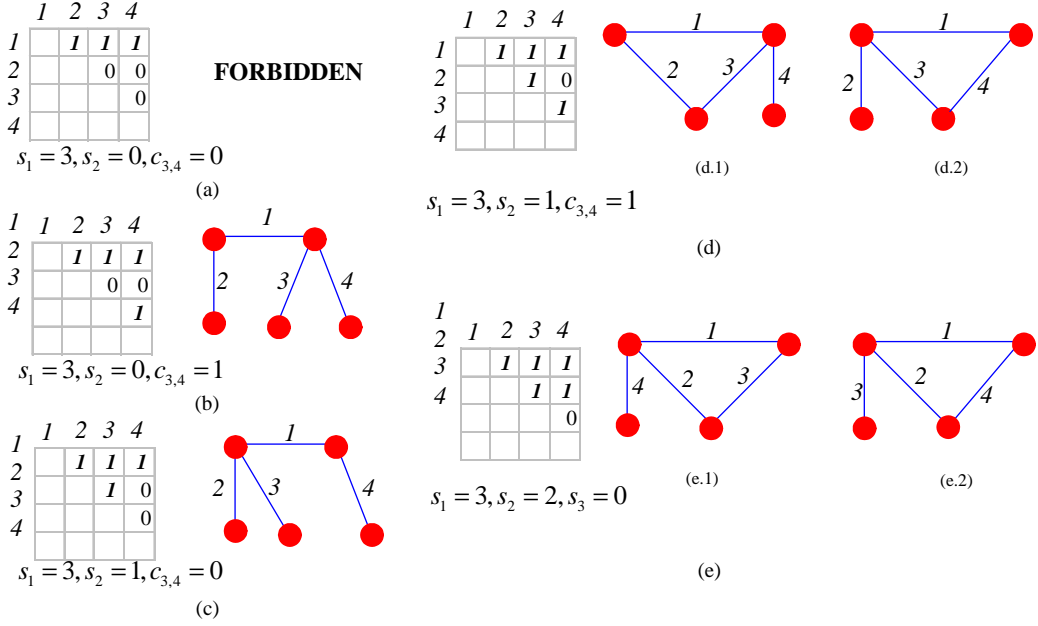}%
\caption{The adjacency patterns of link $1$ and its neighboring links when
$s_{1}=3$. Pattern (a) is forbidden, and patterns (b), (c) and (f) correspond
to only one configuration respectively. Patterns (d) and (e) both have two
possible configurations.}%
\label{s1smaller4_3}%
\end{center}
\end{figure}

There are two recognizable adjacency patterns as described in Figure
\ref{s1smaller4_3} (b), and (c). Taking pattern (c) as an example, links $1$,
$2$ and $3$ are pairwise adjacent, then the configuration of them is $K_{3}$
or $K_{1,3}$, as shown in Figure \ref{s1smaller4_2} (b). Link $4$ is also
adjacent to link $1$, but not adjacent to links $2$ and $3$, suggesting that
the configuration of links $1$, $2$ and $3$ must be $K_{3}$, and link $4$ is
incident to the other endnode of link $1$. Figure \ref{s1smaller4_3} (a)
depicts the smallest forbidden link adjacency pattern in a LAM. The adjacency
relation of links $1$, $2$ and $3$ is recognizable, and the configuration is a
path on four nodes, as shown in Figure \ref{s1smaller4_2} (a). Link $4$ is
adjacent to link $1$, then link $4$ must be also adjacent to links $2$ or $3$.
Hence the pattern is forbidden. If $s_{1}=3$ and $s_{2}=0$, the initial state
is $\mathcal{E}_{3,b}$ (lines 1-2 of Algorithm \ref{algo_initialization_s1_3}%
). If $s_{1}=3$, $s_{2}=1$ and $c_{3,4}=0$, the initial state is
$\mathcal{E}_{3,c}$ (lines 3-5 of Algorithm \ref{algo_initialization_s1_3}).
When $s_{1}=3$, $s_{2}=1$ and $c_{3,4}=1$, due to the connectivity of the
concerned graph, either $c_{2,5}=c_{3,5}=c_{4,5}=1$ or $c_{2,5}=c_{3,5}%
=1,c_{4,5}=0$ or $c_{2,5}=c_{3,5}=0,c_{4,5}=1$ or $c_{2,5}=1,c_{3,5}%
=c_{4,5}=0$ or $c_{2,5}=0,c_{3,5}=c_{4,5}=1$. If $c_{2,5}\neq c_{3,5}$, the
initial state is $\mathcal{E}_{3,d.2}$, else if $c_{2,5}=c_{3,5}\neq c_{4,5}$,
the initial state is $\mathcal{E}_{3,d.1}$, else if $c_{2,5}=c_{3,5}%
=c_{4,5}=1$, we need to look further at the relation of $c_{2,6}$ and
$c_{3,6}$: if $c_{2,6}\neq c_{3,6}$, the initial state is $\mathcal{E}%
_{3,d.2}$, else the initial state is $\mathcal{E}_{3,d.1}$ (lines 11-15 of
Algorithm \ref{algo_initialization_s1_3}). If there are only 5 links in total
and $c_{2,5}=c_{3,5}=c_{4,5}=1$, one can choose any of $\mathcal{E}_{3,d.1}$
and $\mathcal{E}_{3,d.2}$ as the initial state, and get isomorphic
configurations. If $s_{1}=3$, $s_{2}=2$ and $s_{3}=0$, the same method is
employed (lines 21-26 of Algorithm \ref{algo_initialization_s1_3}).%
\begin{equation}
\mathcal{E}_{3,b}=\left[
\begin{array}
[c]{ccccccc}%
1 & 1 & 2 & 2 & 0 & \cdots & 0\\
2 & 0 & 0 & 0 & 0 & \cdots & 0
\end{array}
\right]  ,\mathcal{E}_{3,c}=\left[
\begin{array}
[c]{ccccccc}%
1 & 1 & 1 & 2 & 0 & \cdots & 0\\
2 & 0 & 0 & 0 & 0 & \cdots & 0
\end{array}
\right]
\end{equation}%
\begin{equation}
\mathcal{E}_{3,d.1}=\left[
\begin{array}
[c]{ccccccc}%
1 & 1 & 2 & 2 & 0 & \cdots & 0\\
2 & 0 & 0 & 0 & 0 & \cdots & 0
\end{array}
\right]  ,\mathcal{E}_{3,d.2}=\left[
\begin{array}
[c]{ccccccc}%
1 & 1 & 1 & 2 & 0 & \cdots & 0\\
2 & 0 & 0 & 0 & 0 & \cdots & 0
\end{array}
\right]
\end{equation}%
\begin{equation}
\mathcal{E}_{3,e.1}=\left[
\begin{array}
[c]{ccccccc}%
1 & 1 & 2 & 1 & 0 & \cdots & 0\\
2 & 0 & 0 & 0 & 0 & \cdots & 0
\end{array}
\right]  ,\mathcal{E}_{3,e.2}=\left[
\begin{array}
[c]{ccccccc}%
1 & 1 & 1 & 2 & 0 & \cdots & 0\\
2 & 0 & 0 & 0 & 0 & \cdots & 0
\end{array}
\right]
\end{equation}
%

\begin{algorithm}%

\caption{$E_{2 \times L_{G}} \Leftarrow Initialization3(C,s_{2},s_{3})$}\label{algo_initialization_s1_3}%

\begin{algorithmic}[1]
\IF{$s_{2}=0$}
\STATE $E \Leftarrow \mathcal{E}_{3,b}$
\ELSIF{$s_{2}=1$ and $c_{3,4}=0$}
\STATE $E \Leftarrow \mathcal{E}_{3,c}$
\ELSIF{$s_{2}=1$ and $c_{3,4}=1$}
\IF{$c_{2,5} \neq c_{3,5}$ or ($c_{2,5} = c_{3,5}$ and $c_{2,5} = c_{4,5}$ and $L_{G} = 5$)}
\STATE $E \Leftarrow \mathcal{E}_{3,d.2}$
\ELSIF{$c_{2,5} = c_{3,5}$ and $c_{2,5} \neq c_{4,5}$}
\STATE $E \Leftarrow \mathcal{E}_{3,d.1}$
\ELSIF{$c_{2,5} = c_{3,5}$ and $c_{2,5} = c_{4,5}$ and $c_{2,6} = c_{3,6}$}
\STATE $E \Leftarrow \mathcal{E}_{3,d.1}$
\ELSIF{$c_{2,5} = c_{3,5}$ and $c_{2,5} = c_{4,5}$ and $c_{2,6} \neq c_{3,6}$}
\STATE $E \Leftarrow \mathcal{E}_{3,d.2}$
\ENDIF
\ELSIF{$s_{2}=2$ and $s_{3}=0$}
\IF{$c_{2,5} \neq c_{3,5}$ or ($c_{2,5} = c_{3,5}$ and $c_{2,5} = c_{4,5}$ and $L_{G} = 5$)}
\STATE $E \Leftarrow \mathcal{E}_{3,e.2}$
\ELSIF{$c_{2,5} = c_{3,5}$ and $c_{2,5} \neq c_{4,5}$}
\STATE $E \Leftarrow \mathcal{E}_{3,e.1}$
\ELSIF{$c_{2,5} = c_{3,5}$ and $c_{2,5} = c_{4,5}$ and $c_{2,6} = c_{3,6}$}
\STATE $E \Leftarrow \mathcal{E}_{3,e.1}$
\ELSIF{$c_{2,5} = c_{3,5}$ and $c_{2,5} = c_{4,5}$ and $c_{2,6} \neq c_{3,6}$}
\STATE $E \Leftarrow \mathcal{E}_{3,e.2}$
\ENDIF
\ENDIF
\end{algorithmic}%

\end{algorithm}%

\subsection{When $s_{1}\geq4$}

\subsubsection{When $s_{2}\geq3$}

The configuration is unique. The initial state of $E$ is $\mathcal{E}_{4}$.%
\begin{figure}
[ptb]
\begin{center}
\includegraphics[
height=1.9216in,
width=2.469in
]%
{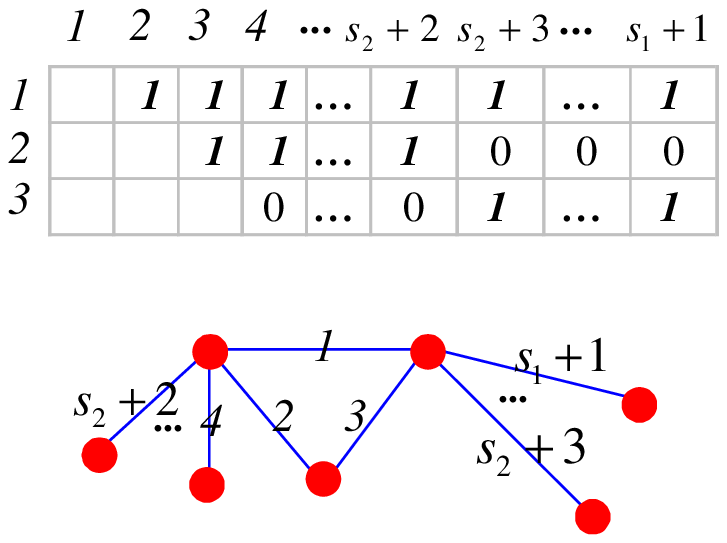}%
\caption{The adjacency pattern and the corresponding configuration when
$s_{3}=0$, $s_{2}\geq3$ and $s_{1}\geq4$.}%
\label{s1_geq4_s2_geq3_s3_eq0}%
\end{center}
\end{figure}
\begin{equation}%
\begin{array}
[c]{c}%
\\
\mathcal{E}_{4}=
\end{array}%
\begin{array}
[c]{c}%
\\
\left[
\begin{array}
[c]{c}%
\\
\\
\end{array}
\right.
\end{array}%
\begin{tabular}
[c]{cccccccccccc}%
${\footnotesize 1}$ & ${\footnotesize 2}$ & ${\footnotesize 3}$ &
${\footnotesize 4}$ & ${\footnotesize \cdots}$ & ${\footnotesize s}%
_{2}{\footnotesize +2}$ & ${\footnotesize s}_{2}{\footnotesize +3}$ &
${\footnotesize \cdots}$ & ${\footnotesize s}_{1}{\footnotesize +1}$ &
${\footnotesize s}_{1}{\footnotesize +2}$ & ${\footnotesize \cdots}$ &
{\footnotesize $L_{G}$}\\\hline
\multicolumn{1}{|c}{$1$} & \multicolumn{1}{|c}{$1$} & \multicolumn{1}{|c}{$2$}
& \multicolumn{1}{|c}{$1$} & \multicolumn{1}{|c}{$\cdots$} &
\multicolumn{1}{|c}{$1$} & \multicolumn{1}{|c}{$2$} &
\multicolumn{1}{|c}{$\cdots$} & \multicolumn{1}{|c}{$2$} &
\multicolumn{1}{|c}{$0$} & \multicolumn{1}{|c}{$\cdots$} &
\multicolumn{1}{|c|}{$0$}\\\hline
\multicolumn{1}{|c}{$2$} & \multicolumn{1}{|c}{$0$} & \multicolumn{1}{|c}{$0$}
& \multicolumn{1}{|c}{$0$} & \multicolumn{1}{|c}{$\cdots$} &
\multicolumn{1}{|c}{$0$} & \multicolumn{1}{|c}{$0$} &
\multicolumn{1}{|c}{$\cdots$} & \multicolumn{1}{|c}{$0$} &
\multicolumn{1}{|c}{$0$} & \multicolumn{1}{|c}{$\cdots$} &
\multicolumn{1}{|c|}{$0$}\\\hline
\end{tabular}%
\begin{array}
[c]{c}%
\\
\left.
\begin{array}
[c]{c}%
\\
\\
\end{array}
\right]
\end{array}
\end{equation}

\subsubsection{When $s_{2}\leq2$}

\begin{figure}
[ptb]
\begin{center}
\includegraphics[
height=6.4342in,
width=4.0577in
]%
{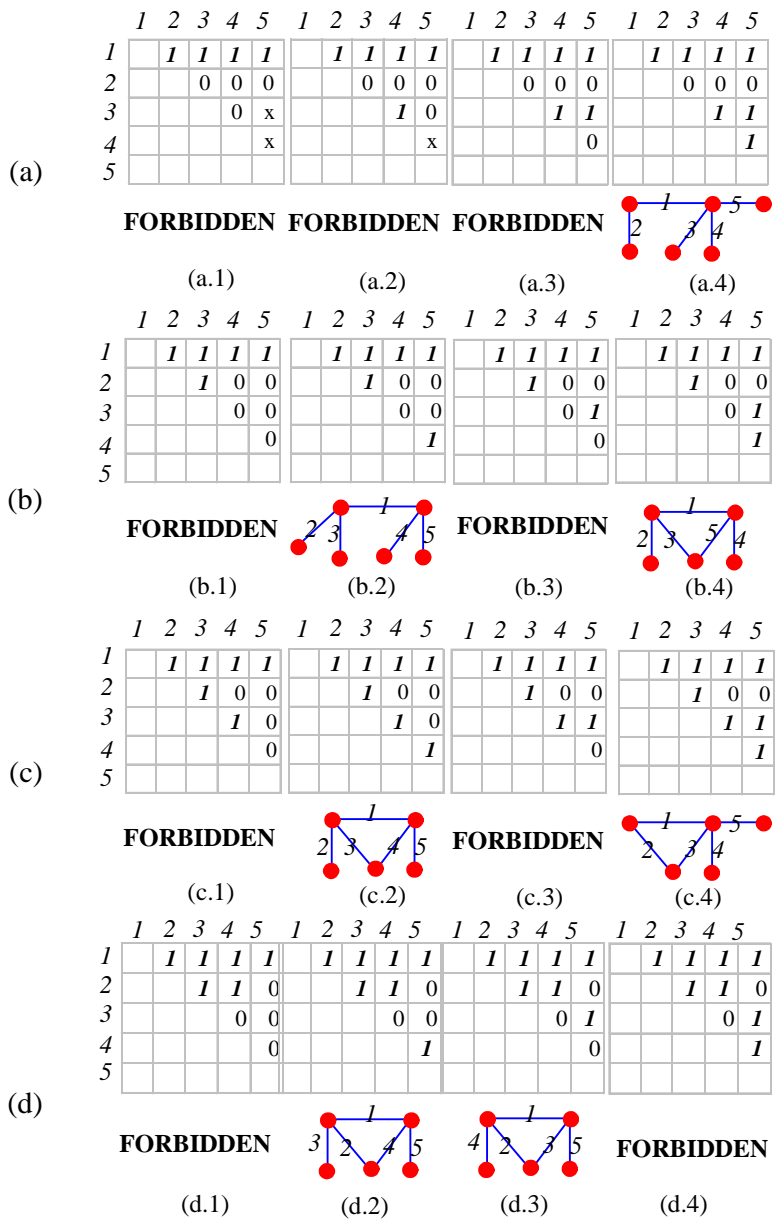}%
\caption{The adjacency patterns of link $1$ and its neighboring links when
$s_{1}=4$. There are $16$ forbidden patterns.\ The other $12$ possible
patterns correspond to only one configuration respectively. The entry x can be
$1$ or $0$.}%
\label{s1smaller4_4}%
\end{center}
\end{figure}

There are $13$ forbidden patterns, as shown in Figure \ref{s1smaller4_4},
where the links with labels larger than $5$ are not displayed. The pattern in
Figure \ref{s1smaller4_3} (a) is forbidden, hence the $4$ patterns in Figure
\ref{s1smaller4_4} (a.1) are also forbidden, where x can be $1$ or $0$. The
pattern of links $1-4$ in Figure \ref{s1smaller4_4} (a.2-3) is the same as the
pattern in Figure \ref{s1smaller4_3} (b), which has a specific configuration.
In Figure \ref{s1smaller4_4} (a.2), link $5$ is adjacent to link $1$ but not
$2$, then link $5$ must be adjacent to link $3$, which is not true, hence the
$2$ patterns in Figure \ref{s1smaller4_4} (a.2) are forbidden. In Figure
\ref{s1smaller4_4} (a.3), link $5$ is adjacent to link $1$ and $3$, then link
$5$ must be adjacent to link $4$, which is not true, hence the pattern in
Figure \ref{s1smaller4_4} (a.3) is also forbidden. Similarly, based on the
patterns in Figure \ref{s1smaller4_3}, we can conclude that patterns in Figure
\ref{s1smaller4_4} (b.1), (b.3), (c.1), (c.3), (d.1) and (d.4) are also
forbidden. Based on the values of entries $s_{2}$, $c_{3,4}$, $c_{3,5}$ and
$c_{4,5}$, Algorithm \ref{algo_initialization_s1_4} decides the initial state
of $E$.%
\begin{equation}
\mathcal{E}_{4,a.4}=\mathcal{E}_{4,c.4}=%
\begin{array}
[c]{c}%
\\
\left[
\begin{array}
[c]{c}%
\\
\\
\end{array}
\right.
\end{array}%
\begin{tabular}
[c]{cccccccccc}%
${\footnotesize 1}$ & ${\footnotesize 2}$ & ${\footnotesize 3}$ &
${\footnotesize 4}$ & ${\footnotesize 5}$ & ${\footnotesize \cdots}$ &
${\footnotesize s}_{1}{\footnotesize +1}$ & ${\footnotesize s}_{1}%
{\footnotesize +2}$ & ${\footnotesize \cdots}$ & {\footnotesize $L_{G}$%
}\\\hline
\multicolumn{1}{|c}{$1$} & \multicolumn{1}{|c}{$1$} & \multicolumn{1}{|c}{$2$}
& \multicolumn{1}{|c}{$2$} & \multicolumn{1}{|c}{$2$} &
\multicolumn{1}{|c}{$\cdots$} & \multicolumn{1}{|c}{$2$} &
\multicolumn{1}{|c}{$0$} & \multicolumn{1}{|c}{$\cdots$} &
\multicolumn{1}{|c|}{$0$}\\\hline
\multicolumn{1}{|c}{$2$} & \multicolumn{1}{|c}{$0$} & \multicolumn{1}{|c}{$0$}
& \multicolumn{1}{|c}{$0$} & \multicolumn{1}{|c}{$0$} &
\multicolumn{1}{|c}{$\cdots$} & \multicolumn{1}{|c}{$0$} &
\multicolumn{1}{|c}{$0$} & \multicolumn{1}{|c}{$\cdots$} &
\multicolumn{1}{|c|}{$0$}\\\hline
\end{tabular}%
\begin{array}
[c]{c}%
\\
\left.
\begin{array}
[c]{c}%
\\
\\
\end{array}
\right]
\end{array}
\end{equation}%
\begin{equation}
\mathcal{E}_{4,b.2}=\mathcal{E}_{4,b.4}=\mathcal{E}_{4,c.2}=\mathcal{E}%
_{4,d.2}=%
\begin{array}
[c]{c}%
\\
\left[
\begin{array}
[c]{c}%
\\
\\
\end{array}
\right.
\end{array}%
\begin{tabular}
[c]{cccccccccc}%
${\footnotesize 1}$ & ${\footnotesize 2}$ & ${\footnotesize 3}$ &
${\footnotesize 4}$ & ${\footnotesize 5}$ & ${\footnotesize \cdots}$ &
${\footnotesize s}_{1}{\footnotesize +1}$ & ${\footnotesize s}_{1}%
{\footnotesize +2}$ & ${\footnotesize \cdots}$ & {\footnotesize $L_{G}$%
}\\\hline
\multicolumn{1}{|c}{$1$} & \multicolumn{1}{|c}{$1$} & \multicolumn{1}{|c}{$1$}
& \multicolumn{1}{|c}{$2$} & \multicolumn{1}{|c}{$2$} &
\multicolumn{1}{|c}{$\cdots$} & \multicolumn{1}{|c}{$2$} &
\multicolumn{1}{|c}{$0$} & \multicolumn{1}{|c}{$\cdots$} &
\multicolumn{1}{|c|}{$0$}\\\hline
\multicolumn{1}{|c}{$2$} & \multicolumn{1}{|c}{$0$} & \multicolumn{1}{|c}{$0$}
& \multicolumn{1}{|c}{$0$} & \multicolumn{1}{|c}{$0$} &
\multicolumn{1}{|c}{$\cdots$} & \multicolumn{1}{|c}{$0$} &
\multicolumn{1}{|c}{$0$} & \multicolumn{1}{|c}{$\cdots$} &
\multicolumn{1}{|c|}{$0$}\\\hline
\end{tabular}%
\begin{array}
[c]{c}%
\\
\left.
\begin{array}
[c]{c}%
\\
\\
\end{array}
\right]
\end{array}
\end{equation}%
\begin{equation}
\mathcal{E}_{4,d.3}=%
\begin{array}
[c]{c}%
\\
\left[
\begin{array}
[c]{c}%
\\
\\
\end{array}
\right.
\end{array}%
\begin{tabular}
[c]{cccccccccc}%
${\footnotesize 1}$ & ${\footnotesize 2}$ & ${\footnotesize 3}$ &
${\footnotesize 4}$ & ${\footnotesize 5}$ & ${\footnotesize \cdots}$ &
${\footnotesize s}_{1}{\footnotesize +1}$ & ${\footnotesize s}_{1}%
{\footnotesize +2}$ & ${\footnotesize \cdots}$ & {\footnotesize $L_{G}$%
}\\\hline
\multicolumn{1}{|c}{$1$} & \multicolumn{1}{|c}{$1$} & \multicolumn{1}{|c}{$2$}
& \multicolumn{1}{|c}{$1$} & \multicolumn{1}{|c}{$2$} &
\multicolumn{1}{|c}{$\cdots$} & \multicolumn{1}{|c}{$2$} &
\multicolumn{1}{|c}{$0$} & \multicolumn{1}{|c}{$\cdots$} &
\multicolumn{1}{|c|}{$0$}\\\hline
\multicolumn{1}{|c}{$2$} & \multicolumn{1}{|c}{$0$} & \multicolumn{1}{|c}{$0$}
& \multicolumn{1}{|c}{$0$} & \multicolumn{1}{|c}{$0$} &
\multicolumn{1}{|c}{$\cdots$} & \multicolumn{1}{|c}{$0$} &
\multicolumn{1}{|c}{$0$} & \multicolumn{1}{|c}{$\cdots$} &
\multicolumn{1}{|c|}{$0$}\\\hline
\end{tabular}%
\begin{array}
[c]{c}%
\\
\left.
\begin{array}
[c]{c}%
\\
\\
\end{array}
\right]
\end{array}
\end{equation}
%

\begin{algorithm}%

\caption{$E_{2 \times L_{G}} \Leftarrow Initialization4(C,s_{1},s_{2},s_{3})$}\label{algo_initialization_s1_4}%

\begin{algorithmic}[1]
\IF{$s_{2} \geq 3$}
\STATE $E \Leftarrow \mathcal{E}_{4}$
\ELSIF{$s_{2}=0$ or ($s_{2}=1$ and $c_{3,4}=1$ and $c_{3,5}=1$ and $c_{4,5}=1$)}
\STATE $E \Leftarrow \mathcal{E}_{4,a.4}$
\ELSIF{$s_{2}=1$ and $c_{3,4}=0$ and $c_{4,5}=1$}
\STATE $E \Leftarrow \mathcal{E}_{4,b.2}$
\ELSIF{$s_{2}=1$ and $c_{3,4}=1$ and $c_{3,5}=0$ and $c_{4,5}=1$}
\STATE $E \Leftarrow \mathcal{E}_{4,c.2}$
\ELSIF{$s_{2}=2$ and $c_{4,5}=1$}
\STATE $E \Leftarrow \mathcal{E}_{4,d.2}$
\ELSIF{$s_{2}=2$ and $c_{4,5}=0$}
\STATE $E \Leftarrow \mathcal{E}_{4,d.3}$
\ENDIF
\end{algorithmic}%

\end{algorithm}%

\end{document}